\documentclass[11pt,a4paper]{article}

\usepackage[T1]{fontenc}
\usepackage[utf8]{inputenc}
\usepackage{mathtools,amsthm,bm}
\usepackage{newtxtext,newtxmath}
\usepackage{microtype}
\usepackage{graphicx}
\usepackage{float}
\usepackage{booktabs,tabularx,array}
\usepackage{enumitem}
\usepackage{geometry}
\usepackage{cancel}
\usepackage[numbers,sort&compress]{natbib}
\usepackage{xcolor}
\definecolor{dgorange}{RGB}{230,105,0}
\definecolor{dggreen}{RGB}{0,125,65}

\definecolor{dgviolet}{RGB}{105,45,150}
\definecolor{dgblue}{RGB}{20,75,160}

\definecolor{dgresp}{RGB}{0,125,65}

\usepackage[colorlinks=true,linkcolor=blue!45!black,citecolor=blue!45!black,urlcolor=blue!45!black]{hyperref}
\usepackage{aliascnt}
\usepackage[nameinlink,capitalise,noabbrev]{cleveref}
\usepackage{subcaption}
\usepackage[normalem]{ulem}
\newtheorem{theorem}{Theorem}[section]
\newaliascnt{proposition}{theorem}
\newtheorem{proposition}[proposition]{Proposition}
\aliascntresetthe{proposition}
\newaliascnt{lemma}{theorem}

\aliascntresetthe{lemma}
\newaliascnt{corollary}{theorem}
\newtheorem{corollary}[corollary]{Corollary}
\aliascntresetthe{corollary}
\newaliascnt{assumption}{theorem}
\newtheorem{assumption}[assumption]{Assumption}
\aliascntresetthe{assumption}
\theoremstyle{definition}
\newaliascnt{definition}{theorem}
\newtheorem{definition}[definition]{Definition}
\aliascntresetthe{definition}
\theoremstyle{remark}
\newaliascnt{remark}{theorem}
\newtheorem{remark}[remark]{Remark}
\aliascntresetthe{remark}

\crefname{theorem}{theorem}{theorems}
\Crefname{theorem}{Theorem}{Theorems}
\crefname{proposition}{proposition}{propositions}
\Crefname{proposition}{Proposition}{Propositions}
\crefname{lemma}{lemma}{lemmas}
\Crefname{lemma}{Lemma}{Lemmas}
\crefname{corollary}{corollary}{corollaries}
\Crefname{corollary}{Corollary}{Corollaries}
\crefname{assumption}{assumption}{assumptions}
\Crefname{assumption}{Assumption}{Assumptions}
\crefname{definition}{definition}{definitions}
\Crefname{definition}{Definition}{Definitions}
\crefname{remark}{remark}{remarks}
\Crefname{remark}{Remark}{Remarks}

\newcommand{\R}{\mathbb{R}}
\newcommand{\Rplus}{\mathbb{R}_{+}}
\newcommand{\Pp}{\mathbb{P}}

\newcommand{\Var}{\operatorname{Var}}
\newcommand{\Cov}{\operatorname{Cov}}
\newcommand{\1}{\mathbf{1}}
\newcommand{\dd}{\,\mathrm{d}}
\newcommand{\weak}{\rightsquigarrow}
\newcommand{\as}{\mathrm{a.s.}}

\newcommand{\cF}{\mathcal{F}}

\newcommand{\cB}{\mathcal{B}}
\newcommand{\cI}{\mathcal{I}}

\newcommand{\Cop}{\overline C}
\newcommand{\Lam}{\Lambda}
\newcommand{\Gam}{\Gamma}
\newcommand{\Psiop}{\Psi}
\newcommand{\ind}{\mathrel{\perp\!\!\!\perp}}

\title{\textbf{Max-Stable Survival Copulas under Dependent Censoring:\\Sharp Identification and Efficient Inference}}
\author{Djibril GUEYE$^{1}$ \footnote{djibril.gueye@quantlabs.fr}, Salima HELALI$^{2}$ \footnote{salima.helali@utc.fr} and Modou WADE$^{3}$ \footnote{modou.wade@cyu.fr}\\[.2cm]
$^{1}$ \small QUANTLABS (subsidiary of the Rainbow Parterns Group) F-92200, Neuilly-sur-Seine, France\\
\\[.2cm]
$^{2}$  Univ de Technologie  Compiègne, 
LMAC,  F-60203  Compiègne Cedex, France\\
\\[.2cm]
$^{3}$ CY  Cergy  Paris Université, THEMA, F-95011 Cergy-Pontoise Cedex, France\\
}
\date{}

\begin{document}
\maketitle

\begin{abstract}

We study dependent censoring within a generalized Cox first-hitting-time.
framework based on compensator reductions. In the genuinely dependent case,
we show that max-stability of the associated survival copula is equivalent to a common operational clock, leading to a nonparametric Marshall Olkin family. We establish sharp identification results under complete and binary observation schemes. Under binary observations, the common clock and event loading are point identified, whereas the censoring
loading and the associated latent survival and dependence structures are only partially identified. We develop nonparametric estimation and inference for
the identifiable components, derive an efficient restricted estimator, and propose an observable test of max-stable compatibility. Simulations and a
real-data application illustrate the practical implications of the identification and efficiency results.
\end{abstract}

\medskip\

\section{Introduction}
Dependent censoring is a fundamental challenge in survival analysis and competing risks, as the event time and the censoring time are only partially observed. In the usual right-censoring scheme, one observes the minimum \(Z=T\wedge C\) together with an indicator recording which component occurs first. Under independent censoring, this observation scheme leads to the classical Kaplan--Meier and Nelson--Aalen estimators and the associated counting-process theory. Under dependent censoring, the same data may not contain enough information to recover the latent margins or the dependence structure without additional restrictions; see, e.g., \citet{Tsiatis1975,ZhengKlein1995,RivestWells2001}.
The nonidentifiability issue is classical. \citet{Tsiatis1975} showed that latent failure-time distributions are generally not identified from a minimum and a cause indicator alone. Sharp bounds compatible with observed subdistributions were derived by \citet{Peterson1976}. This problem is closely related to the partial-identification framework of \citet{Manski2003}, \citet{Tamer2010}, and \citet{ChernozhukovHongTamer2007}. A different line of work restores point identification by imposing a dependence model. Copula-graphic estimators recover survival margins under an assumed copula; see \citet{ZhengKlein1995} and \citet{RivestWells2001}. These approaches show that dependent censoring can be handled when the dependence structure is sufficiently constrained or externally specified. The question studied here is different: within a structural max-stable class, which features of the latent model are identified by the observed censored data alone? \\
\ \\
Common-shock models provide a natural structural framework for dependent event times. The Marshall--Olkin model \citep{MarshallOlkin1967} is a canonical example, in which a common shock can affect several components simultaneously and generates a singular component on the diagonal. Such simultaneous events are particularly relevant when tied failure and censoring times have a substantive interpretation. This construction is closely connected to the theory of extreme-value and max-stable copulas, including the Sibuya-type family and Pickands representations \citep{Sibuya1960,Pickands1981,GudendorfSegers2010,HofertVrins2013}. More recently, generalized Marshall--Olkin constructions have been used to model dependent censoring with positive probability of simultaneous failure and censoring; see \citet{EscobarBachHelali2024} and \citet{Helali2025}. Our approach builds on this common-shock perspective but focuses specifically on what can and cannot be identified from the resulting censored observations.\\
\ \\
The present paper connects this common-shock perspective with generalized Cox first-hitting-time models. In the construction of \citet{GueyeJeanblanc2022}, event times are defined as first passages of increasing stochastic clocks through independent exponential thresholds. The multivariate extension of \citet{GueyeQuintosCox2025} provides the reduction framework used here, allowing several event times to share common and idiosyncratic clocks and yielding joint survival laws in terms of predictable reductions. The dynamic Sibuya copula framework of \citet{GueyeQuintos2025} provides the copula layer used in our max-stability characterization and Pickands representation. Related first-passage and threshold-regression approaches include \citet{LeeWhitmore2006} and \citet{Liu2020CompetingRisks}. In this framework, dependence is induced by the stochastic clocks rather than by the exponential thresholds themselves, through shared factors, jumps, or a common environment.\\
\ \\
What is missing from the existing literature is an identification analysis that
simultaneously accounts for the structural common-shock mechanism and the information lost through observation coarsening. Existing approaches to
dependent censoring typically impose a copula or another dependence structure
to recover latent survival distributions; see, for example,
\citet{ZhengKlein1995,RivestWells2001,Wang2012}. Such approaches can restore
point identification when the dependence structure is specified, but the identifying assumptions concern the latent distribution and are generally not
directly testable from the observed censored data
\citep{Tsiatis1975,Peterson1976}. The present paper takes a different perspective by asking which features of the latent model are identified by the
observed law itself, and how this identification changes when simultaneous event-censoring observations are coarsened into a binary mark. This
perspective is related to the general theory of inference under coarsened data \citep{Hei91}, but here the coarsening has a specific structural
consequence for the identification of the common-shock component
Simultaneous event-censoring observations are central to this distinction. In common-shock models, a tie between the event and censoring times contains information about the shared shock component rather than representing merely
a boundary case. When ties are recorded as a separate category, this information contributes to identification. When they are merged with ordinary events, the corresponding information is lost, leading to a different identification region. This observation motivates the comparison between the
complete three-category and binary observation schemes developed below.\\
\ \\
We focus on the static survival law induced by the predictable compensator reductions. Apart from the product-copula case, we show that max-stability of the survival copula is equivalent to the existence of a common operational clock. The model is therefore nonparametric in calendar time through this clock, while its dependence structure is governed by a finite-dimensional set of parameters. The resulting survival copula is of Marshall--Olkin/Sibuya type and admits an explicit Pickands representation.
The statistical contribution is driven by the observation scheme. In the complete three-category experiment, strict event, strict censoring, and simultaneous event--censoring are recorded separately, preserving the common-shock information and yielding point identification under max-stability. In the binary experiment, simultaneous event--censoring is merged with the event category, resulting in a loss of information that prevents recovery of the common-shock component. We derive the sharp identified set in the binary experiment and show that no statistic based solely on the binary observations can further refine it.\\
\ \\
The inference theory distinguishes between unrestricted and max-stable restricted estimation. In the unrestricted binary model, Nelson--Aalen type estimators arise from the structural representations of the observable submeasures, rather than from an independent-censoring assumption. Under observable max-stable compatibility, the binary law factorizes into the distribution of the observed minimum and a Bernoulli mark, leading to a restricted estimator of the point-identified observable quantities. We establish uniform consistency, functional weak convergence, semiparametric efficiency, and a strict variance reduction relative to the unrestricted estimator when the restriction holds. We also address misspecification through a compatibility test based on the independence of the observed minimum and the binary mark. Rejection rules out any max-stable completion, whereas non-rejection does not uniquely determine the latent completion.
The numerical section is deliberately concise. It validates the multishock data-generating mechanism, illustrates sharp nonidentification, and applies the proposed methodology to the veteran cancer data. The framework distinguishes the identified set under binary observation, point identification under the complete three-category scheme, and sensitivity to the treatment of ties.\\
\ \\
The remainder of the paper is organized as follows. Section~\ref{sec:cox}
introduces the generalized Cox reductions and the two observation schemes. Section~\ref{sec:maxstability} establishes the max-stability characterization and its connection with a common operational clock. Section~\ref{sec:ident}
studies point and set identification under the complete and binary
observation schemes. Section~\ref{sec:inference} develops the corresponding nonparametric estimation and asymptotic theory, including efficiency results.
Section~\ref{sec:testing} proposes a test of observable max-stable
compatibility. Section~\ref{sec:numerical} presents simulation results and an application to the veteran data. Section~\ref{sec:conclusion} concludes
with a discussion of the main implications and limitations of the proposed framework. Proofs are collected in Section~\ref{sec:proofs}.

\section{Generalized Cox Reductions and Observation Schemes}
\label{sec:cox}
\subsection{Generalized Cox construction}
We start from the generalized Cox construction of \citet{GueyeJeanblanc2022} and its multivariate extension in \citet{GueyeQuintosCox2025}. Let $(\Omega,\mathcal A,\mathbb P,\mathbb F)$ be a filtered probability space,
where $\mathbb F=(\mathcal F_t)_{t\ge0}$ is a filtration satisfying the usual
conditions. Let $E_1$ and $E_2$ be two unit-exponential random variables, independent of each other and independent of $\cF_\infty$. Let
$K^1$ and $K^2$ be nondecreasing, cadlag, $\mathbb F$-adapted processes with $K^1_0=K^2_0=0$. The event time and the censoring time are defined by
\begin{equation}\label{eq:first-hitting-TC}
T=\inf\{t\ge0:K^1_t\ge E_1\},
\qquad
C=\inf\{t\ge0:K^2_t\ge E_2\}.
\end{equation}
The exponential thresholds are independent. The times $T$ and $C$, however, need not be independent, because the stochastic clocks may share common factors,
common jumps, or a common environment.
For $J\in\{1,2,12\}$, define the conditional survival processes
\begin{equation}\label{eq:azema-processes}
G^1_t=\Pp(T>t\mid\cF_t),
\qquad
G^2_t=\Pp(C>t\mid\cF_t),
\qquad
G^{12}_t=\Pp(T>t,C>t\mid\cF_t).
\end{equation}
In the generalized Cox framework, these Azema supermartingales admit multiplicative decompositions of the form (see \cite{GueyeJeanblanc2022, GueyeQuintosCox2025}) \
\begin{equation}\label{eq:multiplicative-decomp-main}
G^J_t=\eta^J_t\exp\{-\Lam^J_t\},
\qquad J\in\{1,2,12\}.
\end{equation}
Here $\eta^J$ is a positive martingale normalized at one, and $\Lam^J$ is the predictable reduction associated with the $\mathbb F$-compensator of the
corresponding survival event. In particular, each $\Lam^J$ is predictable, nondecreasing, starts from zero, and diverges at infinity under the usual
properness condition. These predictable reductions are the structural objects carried from the generalized Cox model into the statistical analysis below.

\subsection{Static deterministic reduction model}
The statistical model studied in the rest of the paper is the deterministic
static specialization of the preceding reduction framework. We assume that
$$
\Lam^J_t=\Lam^J(t),\qquad J\in\{1,2,12\},
$$
where the functions are continuous, non decreasing, locally absolutely continuous
on the interior of the support, vanish at zero, and diverge at infinity when the
corresponding lifetime is required to be proper. \\
\ \\
The assumption that the predictable reductions are deterministic is satisfied
in many classical survival and reliability models. The following examples
illustrate typical situations.
\paragraph{Example 1 (Constant hazards).}
Suppose that the event and censoring mechanisms evolve in a homogeneous
environment with constant failure rates. This corresponds, for example, to
electronic components operating under stable laboratory conditions. The
predictable reductions are
\[
\Lambda^1(t)=\lambda_1 t,\qquad
\Lambda^2(t)=\lambda_2 t,\qquad
\Lambda^{12}(t)=\lambda_{12} t,
\]
where $\max(\lambda_1,\lambda_2)\le \lambda_{12}\le \lambda_1+\lambda_2$.

\paragraph{Example 2 (Weibull model).}
In many biomedical and engineering applications, the failure rate changes with
time. A flexible choice is
\[
\Lambda^1(t)=\alpha_1t^\beta,\qquad
\Lambda^2(t)=\alpha_2t^\beta,\qquad
\Lambda^{12}(t)=\alpha_{12}t^\beta,
\]
where $\beta>0$ and
$\max(\alpha_1,\alpha_2)\le \alpha_{12}\le \alpha_1+\alpha_2$.
The corresponding hazards are
\[
(\Lambda^J)'(t)=\alpha_J\beta t^{\beta-1},
\qquad J\in\{1,2,12\}.
\]
For $\beta>1$, the hazard increases with time, as is often observed for ageing
equipment or chronic diseases, whereas $\beta<1$ corresponds to early failures.
\paragraph{Example 3 (Proportional reductions).}
More generally, let $\Lambda_0$ be any cumulative baseline reduction and write
\[
\Lambda^1(t)=\theta_1\Lambda_0(t),\qquad
\Lambda^2(t)=\theta_2\Lambda_0(t),\qquad
\Lambda^{12}(t)=\theta_{12}\Lambda_0(t),
\]
with $\max(\theta_1,\theta_2)\le\theta_{12}\le\theta_1+\theta_2$.
This includes the exponential model
($\Lambda_0(t)=t$), the Weibull model
($\Lambda_0(t)=t^\beta$), and many other parametric families. Such a
representation naturally arises when the event time and censoring time are
affected by the same underlying time scale but with different intensities.\\
\ \\
Since the martingale factors in
\eqref{eq:multiplicative-decomp-main} are normalized, the unconditional
survival functions are
\begin{equation}\label{eq:static-survivals}
S_T(t):=\Pp(T>t)=e^{-\Lam^1(t)},\quad
S_C(t):=\Pp(C>t)=e^{-\Lam^2(t)},\quad
H(t):=\Pp(T>t,C>t)=e^{-\Lam^{12}(t)}.
\end{equation}

The interaction reduction is
\begin{equation}\label{eq:Gamma-def}
\Gam(t)=\Lam^1(t)+\Lam^2(t)-\Lam^{12}(t).
\end{equation}
It measures the departure from additivity of the two marginal reductions on the
diagonal. The induced static bivariate survival law is
\begin{equation}\label{eq:joint-survival}
S_{T,C}(t,s):=\Pp(T>t,C>s)=\exp\{-\Lam^1(t)-\Lam^2(s)+\Gam(t\wedge s)\}.
\end{equation}
This formula is consistent with \eqref{eq:static-survivals}, since setting
$s=t$ gives $\Pp(T>t,C>t)=e^{-\Lam^{12}(t)}$. It contains independent
censoring as the special case $\Gam\equiv0$ and common-shock dependence when
$\Gam$ is nonzero.

\begin{assumption}[Static regularity and admissibility]\label{ass:regular}
The deterministic reductions $\Lam^1,\Lam^2,\Lam^{12}:\Rplus\to\Rplus$ are
continuous, locally absolutely continuous, nondecreasing, and satisfy
$\Lam^J(0)=0$. Define $\Gam=\Lam^1+\Lam^2-\Lam^{12}$. For Lebesgue-a.e.
$t$,
\begin{equation}\label{eq:admissibility}
\Gam'(t)\ge0,
\qquad
(\Lam^1)'(t)-\Gam'(t)\ge0,
\qquad
(\Lam^2)'(t)-\Gam'(t)\ge0.
\end{equation}
\end{assumption}
\noindent
The first set of properties--monotonicity, normalization at zero, and divergence
under properness--comes from the predictable compensator-reduction construction.
The inequalities in \eqref{eq:admissibility} are additional static
admissibility conditions: they ensure that the simultaneous-shock, strict-event,
and strict-censoring submeasures are nonnegative. Equivalently,
\begin{equation}\label{eq:admissibility-equiv}
0\le (\Lam^1)'\le (\Lam^{12})',
\qquad
(\Lam^{12})'-(\Lam^1)'\le (\Lam^2)'\le (\Lam^{12})'
\quad\text{a.e.}
\end{equation}
All identities below have a Lebesgue--Stieltjes version. Absolute continuity is
used only to keep the notation short.

\subsection{Observation schemes and observable submeasures}

Let $Z=T\wedge C$. We distinguish two models. The complete model
records
$$
(Z,\Delta),\qquad
\Delta=\begin{cases}
1,&T<C,\\
0,&C<T,\\
2,&T=C.
\end{cases}
$$
The binary model records
$$
(Z,\delta),\qquad
\delta=\1_{\{T\le C\}}=\1_{\{\Delta=1\}}+\1_{\{\Delta=2\}}.
$$
Thus the binary observation merges strict events and simultaneous events. This
coarsening is the source of partial identification.
The next identities are the basic observed subdistribution identities. They are
written in reduction form because identification is determined by which
submeasures remain visible after coarsening; compare the classical competing-risk
nonidentifiability results of \citet{Tsiatis1975} and \citet{Peterson1976}.

\begin{proposition}[Observable submeasures]\label{prop:submeasures}
Under \Cref{ass:regular}, the density of $Z$ and the complete-category
subdensities are
\begin{align}
 f_Z(t)&=H(t)(\Lam^{12})'(t),\label{eq:fz}\\
 f_{Z,\Delta=1}(t)&=H(t)\{(\Lam^1)'(t)-\Gam'(t)\}=H(t)\{(\Lam^{12})'(t)-(\Lam^2)'(t)\},\label{eq:f1}\\
 f_{Z,\Delta=0}(t)&=H(t)\{(\Lam^2)'(t)-\Gam'(t)\}=H(t)\{(\Lam^{12})'(t)-(\Lam^1)'(t)\},\label{eq:f0}\\
 f_{Z,\Delta=2}(t)&=H(t)\Gam'(t).\label{eq:f2}
\end{align}
Consequently, in the binary model,
\begin{equation}\label{eq:binary-submeasures}
 f_{Z,\delta=1}(t)=H(t)(\Lam^1)'(t),\qquad
 f_{Z,\delta=0}(t)=H(t)\{(\Lam^{12})'(t)-(\Lam^1)'(t)\}
\end{equation}
\end{proposition}

\section{Max-Stability and the Common Operational Clock}
\label{sec:maxstability}

Let $u=e^{-\Lam^1(t)}$ and $v=e^{-\Lam^2(s)}$. When $\Lam^1$ and $\Lam^2$ are strictly increasing and unbounded, define $t_u=(\Lam^1)^{-1}(-\log u)$ and $s_v=(\Lam^2)^{-1}(-\log v)$. The static survival copula is
\begin{equation}\label{eq:surv-copula}
\Cop(u,v)=uv\exp\{\Gam(t_u\wedge s_v)\}.
\end{equation}

\begin{definition}[Max-stability]
A survival copula $\Cop$ is max-stable if $\Cop(u^r,v^r)=\Cop(u,v)^r$ for all $u,v\in(0,1]$ and $r>0$. Equivalently, $\ell(x,y)=-\log\Cop(e^{-x},e^{-y})$ is homogeneous of order one.
\end{definition}
\noindent
The following result is a reduction-based counterpart of the standard max-stability theory of extreme-value copulas \citep{Pickands1981,GudendorfSegers2010,HofertVrins2013}. The point is not that Marshall--Olkin/Sibuya copulas are max-stable, which is classical, but that max-stability forces proportionality of the generalized Cox reductions in the genuinely dependent bivariate case.

\begin{remark}[Product-copula exception]\label{rem:product-exception}
If $\Gam\equiv0$, then $\Cop(u,v)=uv$ is max-stable for arbitrary marginal reductions. Hence the common-clock characterization cannot be stated without excluding the product case: independence is max-stable even when the two marginal clocks are unrelated.
\end{remark}
The next theorem is the bivariate dependent-censoring specialization of the
connected-component max-stability characterization in \citet{GueyeQuintos2025}.
We keep the proof in the paper in order to make the event--censoring notation
self-contained and to isolate the product-copula exception.
\begin{theorem}[Bivariate max-stability characterization]\label{thm:maxstab}
Assume that $\Lam^1$ and $\Lam^2$ are continuous, strictly increasing, and unbounded, and suppose $\Gam\not\equiv0$. Then $\Cop$ is max-stable if and only if there exist constants $\lambda_1,\lambda_2,\lambda_{12}>0$ and a continuous, strictly increasing, unbounded clock $\psi$, with $\psi(0)=0$, such that
\begin{equation}\label{eq:common-clock-main}
\Lam^1=\lambda_1\psi,
\qquad
\Lam^2=\lambda_2\psi,
\qquad
\Lam^{12}=\lambda_{12}\psi.
\end{equation}
The admissibility constraints are
\begin{equation}\label{eq:lambda-constraints}
\max(\lambda_1,\lambda_2)\le\lambda_{12}\le\lambda_1+\lambda_2.
\end{equation}
\end{theorem}

\begin{remark}[Why strict increase is assumed]\label{rem:strict-increase}
Strict increase and unboundedness are imposed in \Cref{thm:maxstab} to work with ordinary inverse marginal reductions in the copula representation. Flat portions can be handled by generalized inverses, but this only adds notation. The identification mechanism is unchanged on the effective support where the reductions increase.
\end{remark}
\noindent
In the genuinely dependent case, normalize
\begin{equation}\label{eq:canonical-normalization}
\Psiop=\Lam^{12},\qquad a=\lambda_1/\lambda_{12},\qquad b=\lambda_2/\lambda_{12}.
\end{equation}
Then
\begin{equation}\label{eq:normalized-model}
\Lam^1=a\Psiop,
\qquad
\Lam^2=b\Psiop,
\qquad
\Lam^{12}=\Psiop,
\qquad
\Gam=(a+b-1)\Psiop,
\end{equation}
where
\begin{equation}\label{eq:Theta}
0<a\le1,
\qquad
0<b\le1,
\qquad
 a+b\ge1.
\end{equation}

\begin{proposition}[Marshall--Olkin/Sibuya form]\label{prop:copula-pickands}
Set $\theta=a+b-1$. The survival copula is
\begin{equation}\label{eq:MO-copula}
\Cop_{a,b}(u,v)=uv\min\{u^{-\theta/a},v^{-\theta/b}\}.
\end{equation}
Its stable-tail dependence function is
\begin{equation}\label{eq:ell-ab}
\ell_{a,b}(x,y)=x+y-\min\{\theta x/a,\theta y/b\}.
\end{equation}
The Pickands function $A_{a,b}(w)=\ell_{a,b}(w,1-w)$ is piecewise affine with breakpoint $w^*(a,b)=a/(a+b)$.
\end{proposition}
\noindent
This is the survival-copula form of the Marshall--Olkin/Sibuya family \citep{MarshallOlkin1967,Sibuya1960,HofertVrins2013}, obtained here from the common-clock reduction representation.

\section{Identification}
\label{sec:ident}

\subsection{Unrestricted binary identification}

The binary submeasures in \eqref{eq:binary-submeasures} immediately show what the binary experiment can and cannot reveal. The result below is a sharp-identification statement in the sense of \citet{Manski2003}, \citet{Tamer2010}, and \citet{ChernozhukovHongTamer2007}: the displayed set is exactly the set of latent reductions compatible with the same observed binary law.

\begin{theorem}[Sharp binary identification without max-stability]\label{thm:unrestricted-binary}
Within the admissible class of \Cref{ass:regular}, two triplets $(\Lam^1,\Lam^2,\Lam^{12})$ and $(\widetilde\Lam^1,\widetilde\Lam^2,\widetilde\Lam^{12})$ generate the same law of $(Z,\delta)$ if and only if
$$
\Lam^{12}=\widetilde\Lam^{12},\qquad \Lam^1=\widetilde\Lam^1.
$$
Given the identified pair $(\Lam^1,\Lam^{12})$, the sharp identified set for $\Lam^2$ is
\begin{equation}\label{eq:I2-unrestricted}
\cI_2(\Lam^1,\Lam^{12})=\bigl\{L:L(0)=0,
\ (\Lam^{12})'-(\Lam^1)'\le L'\le(\Lam^{12})'\text{ a.e.}\bigr\}.
\end{equation}
Equivalently, $0\le\Gam'\le(\Lam^1)'$ a.e.
\end{theorem}

\subsection{Binary max-stable identification}

Under \eqref{eq:normalized-model}, the binary law satisfies
\begin{equation}\label{eq:binary-factorization}
\Pp(Z\in\dd t,\delta=1)=a e^{-\Psiop(t)}\dd\Psiop(t),\qquad
\Pp(Z\in\dd t,\delta=0)=(1-a)e^{-\Psiop(t)}\dd\Psiop(t).
\end{equation}
Hence $Z\ind\delta$ and $\delta\sim\operatorname{Bernoulli}(a)$.

\begin{theorem}[Sharp identification under binary max-stability]\label{thm:binary-MS}
Assume $0<a<1$. In the normalized max-stable model, the binary law point identifies
\begin{equation}\label{eq:Psi-a-id}
\Psiop(t)=-\log\Pp(Z>t),\qquad a=\Pp(\delta=1),
\end{equation}
but does not point identify $b$. The sharp identified set is
\begin{equation}\label{eq:b-sharp}
\cB(a)=[1-a,1].
\end{equation}
All $b\in\cB(a)$ generate the same binary observable law.

\end{theorem}

\begin{remark}[Geometry of the sharp interval]\label{rem:b-geometry}
The lower bound $1-a$ and the upper bound $1$ are exactly the normalized admissibility constraints. Under the common-clock parametrization, the simultaneous-shock component is $\Gam=(a+b-1)\Psiop$. Nonnegativity of this component gives $a+b\ge1$, hence $b\ge1-a$. The strict-censoring rate is nonnegative only if $\dd\Lam^2\le\dd\Lam^{12}$. Since $\Lam^2=b\Psiop$ and $\Lam^{12}=\Psiop$, this gives $b\le1$. The binary mark observes the sum of strict failures and simultaneous shocks, which identifies $a$, but it does not split that sum into its two latent components. Therefore every $b\in[1-a,1]$ is observationally equivalent in the binary experiment.
\end{remark}
\noindent
The identified set propagates to all latent quantities that depend on $b$.

\begin{proposition}[Sharp propagation of the binary identified set]\label{prop:sharp-propagation}
Under binary max-stability, for every $t$,
\begin{equation}\label{eq:latent-bounds}
\Lam^2(t)\in[(1-a)\Psiop(t),\Psiop(t)],
\qquad
S_C(t)\in[e^{-\Psiop(t)},e^{-(1-a)\Psiop(t)}],
\qquad
\Gam(t)\in[0,a\Psiop(t)].
\end{equation}
For every $(u,v)\in[0,1]^2$ and $w\in[0,1]$, the sharp pointwise envelopes are
\begin{equation}\label{eq:copula-envelope-main}
uv\le\Cop_{a,b}(u,v)\le\min\{v,uv^{1-a}\},
\end{equation}
and
\begin{equation}\label{eq:pickands-envelope-main}
1-\min\{w,a(1-w)\}\le A_{a,b}(w)\le1.
\end{equation}
All bounds are sharp. The dependence coefficient formulas and their sharp intervals are collected in \Cref{app:coefficients}.
\end{proposition}

\begin{corollary}
[Sharp identification of Kendall's tau]
\label{prop:kendall-bounds}
Under binary max-stability, Kendall's rank correlation coefficient between
$T$ and $C$ is
\[
\tau_K(T,C)=a+b-1.
\]
Consequently, its sharp identified set is
\[
\mathcal I_{\tau_K}(a)=[0,a].
\]
In particular, although the observable quantities $\Psi$ and $a$ are point
identified, the dependence between $T$ and $C$ is only partially identified
from the binary observation scheme.
\end{corollary}
\begin{corollary}[No observable refinement]\label{cor:no-refinement}
Fix $(\Psiop,a)$. For any $b_1,b_2\in[1-a,1]$, the induced probability measures on the binary sample space are identical. Therefore no estimator, test, or statistic measurable with respect to an i.i.d. sample of $(Z,\delta)$ can consistently distinguish $b_1$ from $b_2$.
\end{corollary}

\subsection{Complete three-category identification}

The complete experiment restores the common-shock information.

\begin{theorem}[Point identification with three categories]\label{thm:three-category}
In the normalized max-stable model, the complete law of $(Z,\Delta)$ satisfies $Z\ind\Delta$ and
\begin{equation}\label{eq:three-probs-main}
\Pp(\Delta=1)=1-b,
\qquad
\Pp(\Delta=0)=1-a,
\qquad
\Pp(\Delta=2)=a+b-1.
\end{equation}
Consequently,
\begin{equation}\label{eq:three-id-main}
\Psiop(t)=-\log\Pp(Z>t),
\qquad
 a=1-\Pp(\Delta=0),
\qquad
 b=1-\Pp(\Delta=1).
\end{equation}
Thus $\Lam^1,\Lam^2,\Lam^{12},\Gam$, the survival copula, and the Pickands function are point identified.
\end{theorem}

\begin{table}[t]
\centering
\small
\caption{Identification under the two recording schemes in the max-stable model.}
\label{tab:id-summary}
\begin{tabularx}{\textwidth}{>{\raggedright\arraybackslash}p{3.7cm}>{\centering\arraybackslash}X>{\centering\arraybackslash}X}
\toprule
Quantity & Binary $(Z,\delta)$ & Complete $(Z,\Delta)$\\
\midrule
Operational clock $\Psiop$ & point identified & point identified\\
Event loading $a$ & point identified & point identified\\
Censoring loading $b$ & $[1-a,1]$ & point identified\\
Censoring margin $S_C$ & set identified & point identified\\
Interaction $\Gam$ & set identified & point identified\\
Copula and Pickands function & set identified & point identified\\
\bottomrule
\end{tabularx}
\end{table}
\section{Inference}
\label{sec:inference}

The identification results of Section~\ref{sec:ident} determine which
structural reductions can be consistently estimated from the binary observation \((Z,\delta)\). Consequently, the estimation strategy differs according to
whether the unrestricted model or the max-stable common-clock restriction is
assumed. The inference procedures developed below therefore estimate only the point-identified components of the model. 
\subsection{Estimation strategy}
In the unrestricted binary model, \Cref{thm:unrestricted-binary} shows that the
observable experiment point identifies the pair
\((\Lambda^1,\Lambda^{12})\), whereas the censoring reduction
\(\Lambda^2\) is only partially identified.
Consequently, inference focuses on the two identified cumulative reductions.
Equivalently, whenever the Radon--Nikodym derivative exists, define the observable local cause ratio
\[
q(t)=\frac{\dd\Lambda^1}{\dd\Lambda^{12}}(t).
\]
In the unrestricted model, \(q(\cdot)\) is an arbitrary measurable function satisfying \(0 \le q(t) \le 1\) a.e., so that no structural restriction is imposed on the evolution of the event mechanism relative to the first-exit mechanism.\\
\ \\
Under the max-stable common-clock restriction, \Cref{thm:binary-MS} shows that the binary experiment point identifies only the observable pair \((\Psi,a)\), where \(\Lambda^1=a\Psi\), \(\Lambda^{12}=\Psi\), or equivalently \(q(t)\equiv a\).
Thus the event reduction is completely determined by the common clock and the
constant loading \(a\). The latent censoring loading \(b\) remains unidentified, and its sharp identified set is
\[
b\in[1-a,1].
\]

\subsection{Unrestricted  Nelson-Aalen estimators}
The estimators introduced in this subsection are structural plug-in estimators
derived directly from the observable submeasure identities established in
\Cref{prop:submeasures}. Although they have the same algebraic form as the
classical Nelson--Aalen estimator, their justification is fundamentally
different. They do not rely on independent censoring, but only on the identification identities
\begin{equation}\label{eq:identities}
\mathbb P(Z\in\dd t)
=
H(t)\dd\Lambda^{12}(t),
\qquad
\mathbb P(Z\in\dd t,\delta=1)
=
H(t)\dd\Lambda^1(t).
\end{equation}
Let \((Z_i,\delta_i)_{i=1}^n\) be i.i.d. copies of \((Z,\delta)\). Define
\[
Y_n(t)=\sum_{i=1}^n \mathbf{1}_{\{Z_i\ge t\}},
\qquad
N_n(t)=\sum_{i=1}^n \mathbf{1}_{\{Z_i\le t\}},
\qquad
N_{1,n}(t)=\sum_{i=1}^n \mathbf{1}_{\{Z_i\le t,\delta_i=1\}}.
\]
Replacing the unknown survival function \(H\) in
\eqref{eq:identities} by its empirical counterpart
\(Y_n/n\) yields the natural plug-in estimators of the two identifiable
reductions:
\begin{equation}\label{eq:NA-estimators}
\widetilde\Lambda_n^{12}(t)
=
\int_{[0,t]}
\frac{\dd N_n(s)}{Y_n(s)},
\qquad
\widetilde\Lambda_n^1(t)
=
\int_{[0,t]}
\frac{\dd N_{1,n}(s)}{Y_n(s)}.
\end{equation}
The notation follows the counting-process formulation of
\citet{AndersenEtAl1993}. However, unlike the classical survival setting,
these estimators are not justified by an independent censoring assumption.
They arise directly from the structural submeasure decomposition (\ref{eq:identities}) of the
binary observation model and therefore remain valid under dependent censoring.\\
\ \\
It is important to emphasize that the binary observation scheme does not allow
the estimation of the censoring reduction \(\Lambda^2\). Indeed, only
\((\Lambda^1,\Lambda^{12})\) are point identified, whereas \(\Lambda^2\) remains
set identified through \(\Cref{thm:unrestricted-binary}\). Consequently, the
survival function of the event time and the survival function of the observed
minimum are identifiable and can be estimated as
\[
\widetilde S_T(t)=\exp\{-\widetilde\Lambda_n^1(t)\},
\qquad
\widetilde H(t)=\exp\{-\widetilde\Lambda_n^{12}(t)\}.
\]
In contrast, the marginal survival function of the censoring time $C$ and the full
joint survival distribution of \((T,C)\) cannot be recovered without additional
structural assumptions.

\subsection{Restricted estimation under max-stability}

We now consider inference under the max-stable common-clock restriction. Recall that, under the normalized representation, \(\Lambda^1=a\Psi\) and \(\Lambda^{12}=\Psi\), where \(\Psi\) is the common clock and \(a\in(0,1]\) is the observable event loading. By \Cref{thm:binary-MS}, the binary observation law point identifies the pair \((\Psi,a)\), with
\[
\Psi(t)=-\log\mathbb P(Z>t),
\qquad
a=\mathbb P(\delta=1).
\]
In contrast, the censoring loading \(b\) remains unidentified.\\
\ \\
Since \(\Psi=\Lambda^{12}\), we estimate the common clock using the structural
Nelson--Aalen estimator introduced previously:
\begin{equation}\label{eq:Psi-estimator}
\widetilde \Psi_n(t)
=
\widetilde\Lambda_n^{12}(t)
=
\int_{[0,t]}
\frac{\dd N_n(s)}{Y_n(s)}.
\end{equation}
The event loading is estimated by
\begin{equation}\label{eq:a-estimator}
\widetilde a_n
=
\frac1n\sum_{i=1}^n\delta_i .
\end{equation}
Therefore, the restricted estimator of the event reduction is
\begin{equation}\label{eq:restricted-Lambda1-estimator}
\widetilde\Lambda_n^{1,R}(t)
=
\widehat a_n\,\widehat\Psi_n(t).
\end{equation}
Compared with the unrestricted estimator \(\widetilde\Lambda_n^1\), this estimator imposes the structural restriction $\displaystyle{\dd\Lambda^1/\dd\Lambda^{12}(t)=a}$ instead of allowing a time-varying ratio. Hence, the max-stable restriction pools information over time through the common-clock representation.\\
\ \\
The restriction identifies the event-time survival function, whose estimator is
\[
\widetilde S_T^{\,R}(t)
=
\exp\{-\widetilde a_n\widetilde\Psi_n(t)\}.
\]
However, \(b\) remains unidentified, so the censoring survival function and the joint survival distribution are only set identified. For each admissible value \(b\in[1-\widetilde a_n,1]\),
the corresponding estimators are
\[
\widetilde S_C^{(b)}(t)
=
\exp\{-b\widetilde\Psi_n(t)\},
\]
and
\[
\widetilde S_{T,C}^{(b)}(t,s)
=
\exp\{-\widetilde a_n\widetilde\Psi_n(t)
-b\widetilde\Psi_n(s)
+(\widetilde a_n+b-1)\widetilde\Psi_n(t\wedge s)\}.
\]
Thus, the max-stable restriction yields point estimation of the observable
quantities \((\Psi,a)\), while the latent censoring loading \(b\) remains only
partially identified. The resulting family of estimators provides a sensitivity
analysis over \(b\in[1-\widetilde a_n,1]\), rather than an estimator of \(b\)
itself.\\
\ \\
The following result establishes the asymptotic properties of the restricted
estimators introduced above.
\begin{theorem}[Consistency and functional CLT]\label{thm:FCLT}
Under binary max-stability and $0<a<1$,
$$
\widetilde a_n\to a\quad\as,
\qquad
\sup_{t\le\tau}|\widetilde\Psiop_n(t)-\Psiop(t)|\to0\quad\as,
\qquad
\sup_{t\le\tau}|\widetilde\Lam^{1,R}_n(t)-\Lam^1(t)|\to0\quad\as.
$$
Moreover,
\begin{equation}\label{eq:FCLT}
\left(\sqrt n(\widetilde\Psiop_n-\Psiop),\sqrt n(\widetilde a_n-a)\right)
\weak(\mathbb G_\Psi,G_a)
\end{equation}
in $\ell^\infty[0,\tau]\times\R$, where $G_a\sim N\{0,a(1-a)\}$, independent of $\mathbb G_\Psi$, and
\begin{equation}\label{eq:GPsi-cov}
\Cov\{\mathbb G_\Psi(s),\mathbb G_\Psi(t)\}=\frac{H(s\vee t)-H(s)H(t)}{H(s)H(t)}.
\end{equation}
For fixed $t$,
\begin{equation}\label{eq:VR-main}
\sqrt n\{\widetilde\Lam^{1,R}_n(t)-\Lam^1(t)\}
\weak N\{0,V_R(t)\},
\qquad
V_R(t)=a(1-a)\Psiop(t)^2+a^2\frac{1-H(t)}{H(t)}.
\end{equation}
\end{theorem}

\begin{theorem}[Efficiency and variance gain]\label{thm:efficiency-main}
In the observable product model $P_Z\otimes\operatorname{Bernoulli}(a)$, the efficient influence functions for $a$, $\Psiop(t)$, and $\Lam^1(t)=a\Psiop(t)$ are
\begin{align}
\phi_a(z,d)&=d-a,\label{eq:IF-a-main}\\
\phi_{\Psi,t}(z,d)&=-\{\1_{\{z>t\}}-H(t)\}/H(t),\label{eq:IF-Psi-main}\\
\phi_{1,t}(z,d)&=\Psiop(t)(d-a)-a\{\1_{\{z>t\}}-H(t)\}/H(t).\label{eq:IF-L1-main}
\end{align}
Thus $\widehat\Lam^{1,R}_n(t)$ is semiparametrically efficient for the point-identified observable target. If the max-stable restriction is correct, the unrestricted estimator $\widetilde\Lam^1_n(t)$ has limiting variance
\begin{equation}\label{eq:VU-main}
V_U(t)=a\{e^{\Psiop(t)}-1\}=a\frac{1-H(t)}{H(t)},
\end{equation}
and
\begin{equation}\label{eq:variance-gain-main}
V_U(t)-V_R(t)=a(1-a)\{e^{\Psiop(t)}-1-\Psiop(t)^2\}\ge0,
\end{equation}
with strict inequality for $0<a<1$ and $\Psiop(t)>0$.
\end{theorem}

\subsection{Inference for the identified set}
The previous sections show that, under the binary observation scheme and the
max-stable restriction, the parameters \(\Psi\) and \(a\) are point identified,
whereas the censoring loading \(b\) is only partially identified. More
precisely, the data do not determine a unique value of \(b\), but only the sharp
identified set
\[
\mathcal B(a)=[1-a,1].
\]
Therefore, classical confidence intervals for a point parameter are not
appropriate for \(b\). The relevant inferential objective is instead to
construct a confidence region that contains the whole identified set with a
prescribed asymptotic probability. This is the standard approach in partially
identified models \citep{ChernozhukovHongTamer2007}.\\
\ \\
Let \(U_{a,n}(1-\alpha)\) be a one-sided upper confidence bound for the
identified parameter \(a\), satisfying
\[
\liminf_{n\to\infty}
\mathbb P\{a\le U_{a,n}(1-\alpha)\}\ge1-\alpha .
\]
Since the lower endpoint of \(\mathcal B(a)\) is \(1-a\), an upper confidence
bound for \(a\) directly provides a lower confidence bound for \(1-a\). Hence,
we define
\begin{equation}\label{eq:confidence-set-b}
\mathcal C_{b,n}(1-\alpha)
=
[1-U_{a,n}(1-\alpha),1].
\end{equation}
By construction,
\[
\liminf_{n\to\infty}
\mathbb P\{\mathcal B(a)\subseteq
\mathcal C_{b,n}(1-\alpha)\}
\ge1-\alpha .
\]
Thus, \(\mathcal C_{b,n}(1-\alpha)\) is not a confidence interval for an
unknown value of \(b\), but a confidence region for all values of \(b\) that are
compatible with the observed binary data.\\
\ \\
For interior values of \(a\), the upper confidence bound can be obtained from
the asymptotic normal approximation of \(\widetilde a_n\):
\[
U_{a,n}(1-\alpha)
=
\min\left\{
1,\,
\widetilde a_n+
z_{1-\alpha}
\sqrt{\frac{\widetilde a_n(1-\widetilde a_n)}{n}}
\right\}.
\]
Near the boundary of the parameter space, an exact binomial confidence bound may be preferred.\\
\ \\
The same principle applies to functional quantities. For example, although the censoring reduction \(\Lambda^2(t)=b\Psi(t)\) is not point identified, the model implies the sharp bounds \((1-a)\Psi(t)\le \Lambda^2(t)\le\Psi(t)\).\\
\ \\
Let \(\tau<\infty\) be a fixed study horizon such that \(H(\tau)=\mathbb P(Z>\tau)>0\). Let \([L_{\Psi,n}(t),U_{\Psi,n}(t)]\) be a simultaneous confidence band for \(\Psi\) on \([0,\tau]\). Combining this
band with the confidence bound for \(a\) gives the outer confidence band
\begin{equation}\label{eq:outer-L2-main}
[(1-U_{a,n})L_{\Psi,n}(t),\,U_{\Psi,n}(t)],
\end{equation}
which contains the identified set of \(\Lambda^2(t)\) with asymptotic coverage
at least \(1-\alpha\).\\
\ \\
Similarly, the dependence component satisfies \(0\le \Gamma(t)\le a\Psi(t)\), and therefore admits the outer confidence band
\begin{equation}\label{eq:outer-Gamma-main}
[0,U_{a,n}U_{\Psi,n}(t)].
\end{equation}
The same construction can be applied to other latent quantities, such as the
joint survival distribution or the associated copula functionals. Hence, even
when some model components are not point identified, the max-stable structure allows valid inference on the entire collection of distributions compatible
with the observed binary data.\\
\ \\
The same monotone construction can be applied to other latent functionals of
the model. In particular, it yields outer confidence envelopes for the compatible copula and Pickands dependence functions. A multiplier bootstrap procedure for constructing the simultaneous confidence band
\([L_{\Psi,n},U_{\Psi,n}]\), jointly with the confidence bound
\([0,U_{a,n}]\), is provided in \Cref{app:bootstrap}.

\section{Testing Observable Max-Stable Compatibility}
\label{sec:testing}

The restricted estimation procedure developed in the previous section relies on the observable proportionality restriction \(\Lambda^1=a\Lambda^{12}\),
which is the binary implication of the max-stable common-clock model. Before
using this restricted estimator, it is therefore important to assess whether
this restriction is compatible with the observed binary data.\\
\ \\
Since only the binary observation $(Z,\delta)$
is available, the latent max-stable representation cannot be tested directly.
The data can only test whether the observed binary law admits at least one max-stable latent completion. This is an observable compatibility problem,
rather than a test of uniqueness of the latent model.

\subsection{Observable compatibility characterization}
Recall the local observable cause ratio introduced in the estimation section,
\[
q(t)=\frac{\dd\Lambda^1}{\dd\Lambda^{12}}(t).
\]
By the binary submeasure identities,
\[
\mathbb P(Z\in\dd t,\delta=1)=H(t)\dd\Lambda^1(t),
\qquad
\mathbb P(Z\in\dd t)=H(t)\dd\Lambda^{12}(t),
\]
so that
\[
q(t)=\mathbb P(\delta=1\mid Z=t)
\]
is a version of the conditional probability of observing an event at a given value of the observed time. Under the unrestricted binary model, \(q(t)\) is allowed to vary with time. In contrast, the max-stable common-clock model imposes
\[
\Lambda^1=a\Lambda^{12},
\qquad
q(t)\equiv a.
\]
The next result shows that this observable restriction admits several equivalent characterizations.
\begin{theorem}[Compatibility characterization]\label{thm:compatibility-main}
For an admissible binary law, the following statements are equivalent:
\begin{enumerate}
\item[i)] There exists \(a\in[0,1]\) such that \(\Lambda^1=a\Lambda^{12}\).

\item[ii)] The local observable cause ratio is constant, \(q(t)=a,\ \dd\Lambda^{12}\text{-a.e.}\).

\item[iii)] The observed time and the binary mark are independent, \(Z\ind\delta\).

\item[iv)] The binary law admits at least one max-stable latent completion.
\end{enumerate}
\end{theorem}
\noindent
Thus, testing observable max-stable compatibility reduces to testing the
independence relation between the observed time \(Z\) and the indicator
\(\delta\). Importantly, rejection of this hypothesis rules out the
max-stable restricted model, whereas acceptance only establishes observable
compatibility and does not imply uniqueness of the latent completion.

\subsection{goodness-of-fit test for the compatibility restriction}

The previous theorem suggests a natural empirical test. Let
\[
F(t)=\mathbb P(Z\le t),
\qquad
F_1(t)=\mathbb P(Z\le t,\delta=1).
\]
Under the compatibility restriction, \(F_1(t)=aF(t)\), where \(a=\mathbb P(\delta=1)\).
For a sample \((Z_i,\delta_i)_{i=1}^n\), define the empirical quantities
\[
\widetilde F_n(t)=\frac1n\sum_{i=1}^n\mathbf 1_{\{Z_i\le t\}},
\qquad
\widetilde F_{1,n}(t)=\frac1n\sum_{i=1}^n\delta_i\mathbf 1_{\{Z_i\le t\}},
\qquad
\widetilde a_n=\frac1n\sum_{i=1}^n\delta_i .
\]
The empirical deviation from the compatibility restriction is measured by
\[
\mathbb T_n(t)
=
\sqrt n
\left\{
\widetilde F_{1,n}(t)
-
\widetilde a_n\widetilde F_n(t)
\right\}.
\]

\begin{theorem}[Null limit and consistency]\label{thm:test-main}
Assume that $0<a<1$, that $F$ is continuous, and that
$Z\ind\delta$. Let $\tau<\infty$ be fixed such that $H(\tau)>0$.
Then
\[
\mathbb T_n
\weak
\sqrt{a(1-a)}\,B\circ F
\]
in $\ell^\infty[0,\tau]$, where $B$ is a standard Brownian bridge.
Consequently,
\[
\frac{\sup_{0\le t\le\tau}|\mathbb T_n(t)|}
{\sqrt{\widehat a_n(1-\widehat a_n)}}
\weak
\sup_{0\le t\le\tau}|B(F(t))|.
\]
The resulting Kolmogorov--Smirnov-type test is consistent against fixed
alternatives for which $Z$ and $\delta$ are dependent. Moreover,
conditional on the observed times and on the total number of ones,
permutation of the binary marks provides an exact finite-sample
conditional calibration under the null hypothesis.
\end{theorem}

\subsection{Consequences of misspecification}
The restricted estimator developed previously is valid when the compatibility
restriction holds. If this restriction is violated, it does not converge to the
unrestricted observable reduction \(\Lambda^1\). Instead, it converges to the
proportional approximation imposed by the max-stable model.
Let
\[
a_\star=\mathbb P(\delta=1)
=
\int q(t)e^{-\Psi(t)}\dd\Psi(t).
\]
Then,
\[
\widetilde\Lambda^{1,R}_n(t)
\longrightarrow
a_\star\Psi(t),
\]
and the difference with the true observable reduction is
\[
a_\star\Psi(t)-\Lambda^1(t)
=
\int_{[0,t]}
\{a_\star-q(s)\}\dd\Psi(s).
\]
Hence, when the compatibility condition fails, the restricted estimator should
be interpreted as a projection onto the max-stable common-clock family rather
than as a consistent estimator of the unrestricted reduction.

The compatibility test therefore plays a complementary role to the efficiency
results obtained previously: it verifies the observable restriction that
justifies the use of the efficient restricted estimator, while preserving the
distinction between observable compatibility and latent identification.

\section{Simulation and Real Data Analysis}
\label{sec:numerical}
\subsection{Simulation study}\label{sub1}

The numerical experiments are designed to illustrate the main theoretical findings of the paper. In particular, we aim to validate the identification
results, compare the unrestricted and max-stable restricted estimation strategies, and assess the practical gain obtained by exploiting the
max-stable common-clock structure.\\
\ \\
Throughout the simulations, the latent variables $(T,C)$ are used  for
the generation of the data. The statistical inference is performed from the observed binary data
\[
(Z,\delta),\qquad Z=T\wedge C,\quad \delta=\mathbf{1}_{\{T\le C\}},
\]
which corresponds to the observation scheme considered in the theoretical analysis.
\subsubsection{Unrestricted estimation}
We first consider the unrestricted binary model, without imposing the
max-stable common-clock restriction. According to
\Cref{thm:unrestricted-binary}, the binary observation law identifies only the
two cumulative reductions
\[
\Lambda^1
\quad\text{and}\quad
\Lambda^{12},
\]
whereas the censoring reduction \(\Lambda^2\) is not point identified. To illustrate this partial identification result, we consider a model with a
time-varying local cause ratio
$\displaystyle{q(t)=\frac{\dd\Lambda^1}{\dd\Lambda^{12}}(t)}$.
We choose $\displaystyle{\Lambda^{12}(t)=t^{1.5}}$,
and $\displaystyle{q(t)=0.3+0.4\frac{t}{1+t}}$.
The first identifiable reduction is then defined by
\[
\Lambda^1(t)
=
\int_0^t q(s)\,\dd\Lambda^{12}(s).
\]
The simulated sample consists of $n=100, 500 \text{ and }1000$ independent copies
$(Z_i,\delta_i)_{i=1}^{n}$ of the observed binary experiment, where
$Z=T\wedge C$ and $\delta=\mathbf{1}_{\{T\le C\}}$.
The variable \(Z\) is generated from its survival function
\[
\Pp(Z>t)=\exp\{-\Lambda^{12}(t)\},
\]
and the mark is generated conditionally on \(Z\) according to
\[
\Pp(\delta=1\mid Z=t)=q(t).
\]
We estimate the two identifiable quantities using the structural
Nelson--Aalen estimators
\[
\widetilde\Lambda_n^{12}(t)
=
\int_{[0,t]}\frac{\dd N_n(s)}{Y_n(s)},
\qquad
\widetilde\Lambda_n^{1}(t)
=
\int_{[0,t]}\frac{\dd N_{1,n}(s)}{Y_n(s)}.
\]
The simulation evaluates the convergence of these estimators by comparing their estimated paths with the true reductions and by computing integrated
squared errors for increasing sample sizes.

\begin{figure}[H]
\centering
\begin{minipage}[t]{.48\textwidth}
\centering
\includegraphics[width=\textwidth]{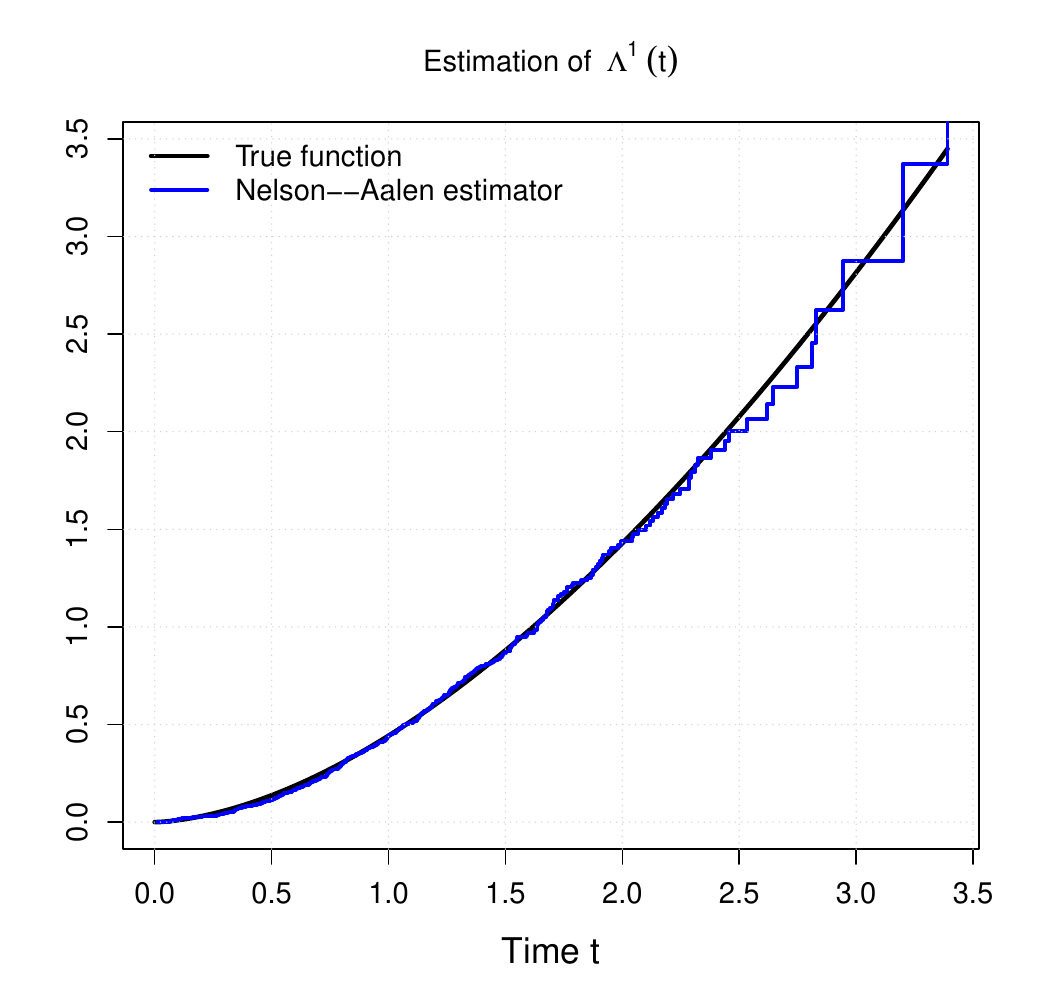}
\end{minipage}\hfill
\begin{minipage}[t]{.48\textwidth}
\centering
\includegraphics[width=\textwidth]{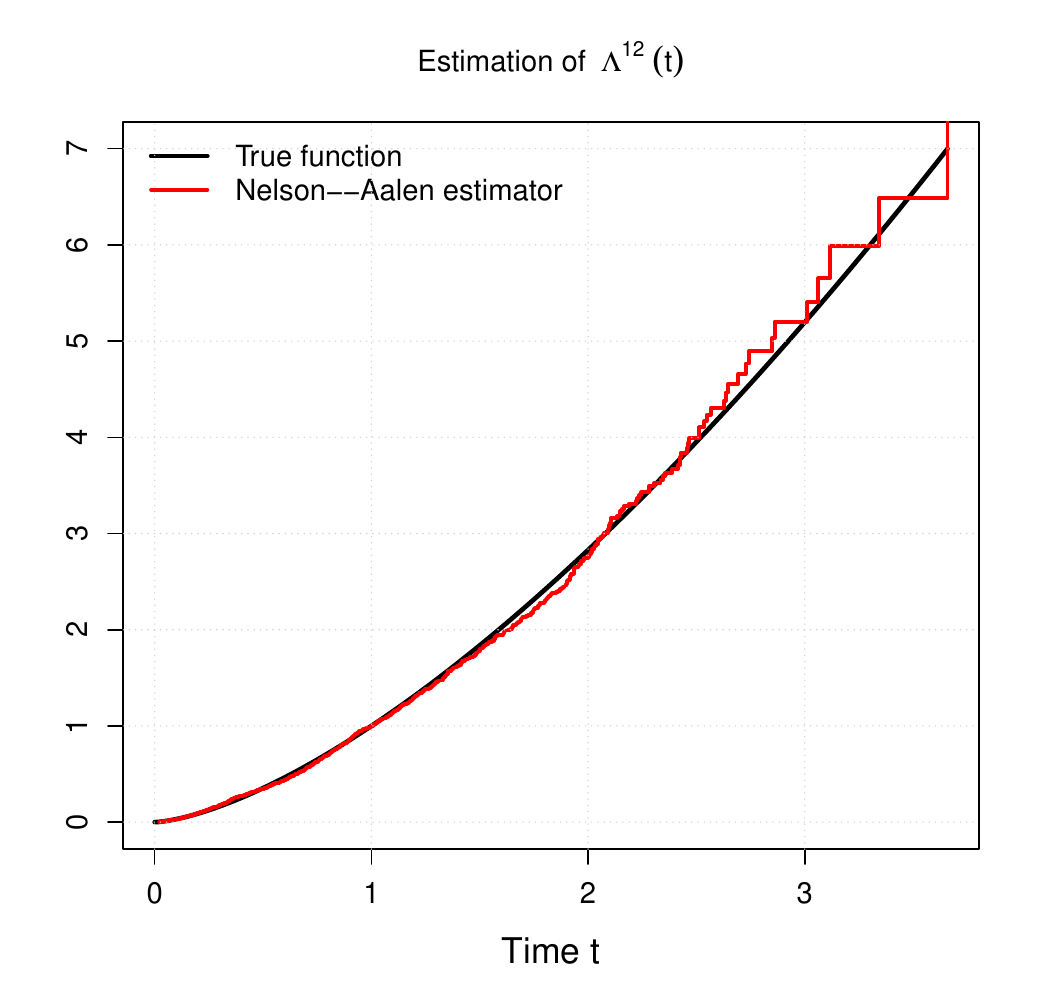}
\end{minipage}
\caption{}
\label{fig:unrestricted-estimation}
\end{figure}

\begin{table}[H]
\centering
\caption{Integrated squared errors (ISE) of the structural Nelson--Aalen estimators under the unrestricted model.}
\label{tab:ISE-unrestricted}
\begin{tabular}{ccc}
\hline
Sample size $n$ &
ISE$(\widetilde{\Lambda}_n^{12})$ &
ISE$(\widetilde{\Lambda}_n^{1})$ \\
\hline
100  & $ 0.074464$ & $0.060053$ \\
500  & $0.040228$ & $0.040228$  \\
1000 & $ \mathbf{0.026168}$ & $  \mathbf{0.021833}$ \\
\hline
\end{tabular}
\end{table}
Figure~\ref{fig:unrestricted-estimation} compares the true cumulative reductions with their structural Nelson--Aalen estimators obtained from the simulated binary observations. In both cases, the estimated curves closely follow the theoretical functions. As the sample size increases, the empirical trajectories become smoother and approach the true reductions uniformly, illustrating the consistency of the structural estimators.
Table~\ref{tab:ISE-unrestricted} reports the integrated squared errors (ISE) for different sample sizes. The estimation error decreases steadily with the sample size for both $\Lambda^{12}$ and $\Lambda^{1}$, confirming the convergence predicted by the theoretical results. As expected from \Cref{thm:unrestricted-binary}, only the two identifiable reductions are estimated in the unrestricted model, whereas no estimator of the nonidentified censoring reduction $\Lambda^{2}$ is considered.

\subsubsection{ Max-stable restricted estimation}
We now investigate the restricted estimation procedure under the max-stable common-clock model. In this setting,
\[
\Lambda^1=a\Psi,
\qquad
\Lambda^{12}=\Psi,
\]
so that the binary experiment point identifies the common clock
\(\Psi\) and the event loading \(a\), whereas the censoring loading \(b\) remains unidentified (\Cref{thm:binary-MS}).\\
\ \\
The objective of this simulation study is threefold:
(i) to evaluate the estimation of the point-identified quantities \(\Psi\), \(a\) and \(\Lambda^1\);
(ii) to illustrate the non-identifiability of \(b\);
and (iii) to compare the restricted estimator with its unrestricted counterpart.\\
\ \\
Two max-stable models are considered.We consider two max-stable common-clock specifications:
\begin{enumerate}
\item[$a)$] The first model corresponds to an exponential operational clock, 
\[
\Psi(t)=2\log(1+t), \qquad a=0.7, \qquad b=0.5.
\]
The associated cumulative reductions are
\[
\Lambda^{12}(t)=2\log(1+t),
\qquad
\Lambda^1(t)=1.4\log(1+t),
\qquad
\Lambda^2(t)=\log(1+t).
\]
The event-time, censoring-time, and joint survival functions are therefore
\[
S_T(t)=\exp\{-1.4\log(1+t)\}=(1+t)^{-1.4},
\qquad
S_C(t)=\exp\{-\log(1+t)\}=(1+t)^{-1},
\]
\[
\Gamma(t)=(a+b-1)\Psi(t)
=0.2\times 2\log(1+t)
=0.4\log(1+t),
\]
and
\[
S_{T,C}(t,s)
=
\exp\left\{
-1.4\log(1+t)
-\log(1+s)
+0.4\log(1+t\wedge s)
\right\}.
\]

\item[$b)$] The second model considers a nonlinear operational clock,
\[
\Psi(t)=t^{3/2}, \qquad a=0.7, \qquad b=0.5.
\]
The associated cumulative reductions are
\[
\Lambda^{12}(t)=t^{3/2},
\qquad
\Lambda^1(t)=0.7t^{3/2},
\qquad
\Lambda^2(t)=0.5t^{3/2}.
\]
The corresponding event-time, censoring-time, and joint survival functions are
\[
S_T(t)=\exp(-0.7t^{3/2}),
\qquad
S_C(t)=\exp(-0.5t^{3/2}),
\]
and
\[
S_{T,C}(t,s)
=
\exp\{-0.7t^{3/2}-0.5s^{3/2}
+0.2(t\wedge s)^{3/2}\}.
\]
\end{enumerate}
For each model, the simulated samples are generated from the multishock
representation implied by the max-stable common-clock restriction. The common shock component is given by
\[
\Gamma
=
\Lambda^1+\Lambda^2-\Lambda^{12}
=
(a+b-1)\Psi .
\]
Therefore, on the operational scale $u=\Psi(t)$,
the three independent shock mechanisms have constant rates
\[
1-b,\qquad 1-a,\qquad a+b-1,
\]
corresponding respectively to an event-only shock, a censoring-only shock,
and a common shock. This representation is the standard multishock
construction associated with the Marshall--Olkin/Sibuya form established in
\Cref{prop:copula-pickands}.
Accordingly, we generate three independent exponential variables
\[
X_1\sim\operatorname{Exp}(1-b),\qquad
X_2\sim\operatorname{Exp}(1-a),\qquad
X_{12}\sim\operatorname{Exp}(a+b-1),
\]
where the rates correspond to the three operational shock intensities.
The operational event and censoring times are then defined as
\[
\tau_T=\min(X_1,X_{12}),
\qquad
\tau_C=\min(X_2,X_{12}).
\]
Since these variables are generated on the common-clock scale, the  times are obtained by applying the inverse clock transformation:
\[
T=\Psi^{-1}(\tau_T),
\qquad
C=\Psi^{-1}(\tau_C).
\]
For each model, independent copies of the observable binary experiment
\[
(Z_i,\delta_i)_{i=1}^{n},
\qquad
Z_i=T_i\wedge C_i,
\qquad
\delta_i=\mathbf1_{\{T_i\le C_i\}},
\]
are generated for sample sizes \(n\in\{50,100,500\}\),
and inference is performed solely from these observable variables.\\
\ \\
For each simulated sample, we estimate the common clock \(\Psi\), the event
loading \(a\), the event reduction \(\Lambda^1=a\Psi\), and the corresponding
event-time survival function. We then investigate the family of compatible
censoring survival functions associated with the identified set
\(b\in[1-a,1]\), and finally compare the estimation accuracy of the restricted
and unrestricted procedures through Monte Carlo experiments based on \(M=1000\) independent replications.

\begin{figure}[H]
\centering

\begin{subfigure}{0.48\textwidth}
    \centering
    \includegraphics[width=\linewidth]{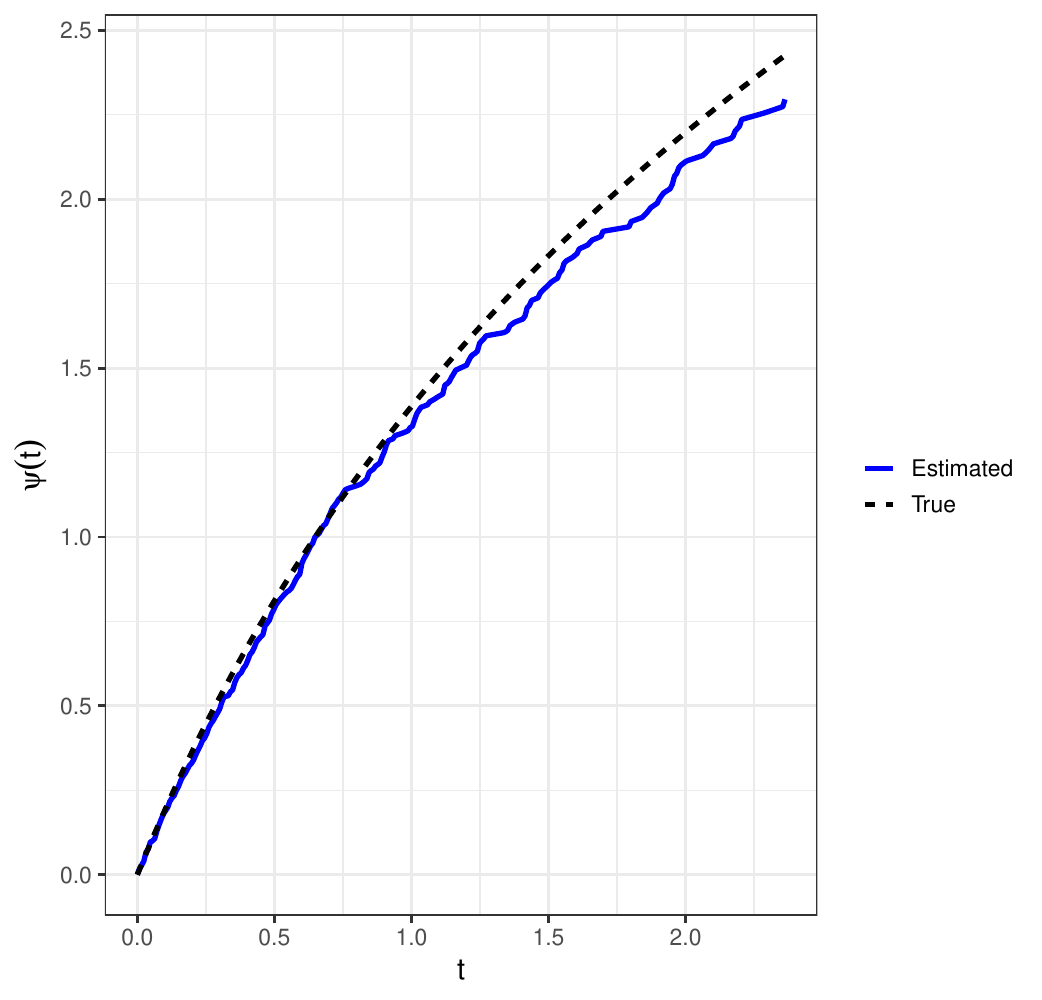}
    \caption{Results for the estimator of $\psi$ for model a).}
    \label{fig:psi-a}
\end{subfigure}
\hfill
\begin{subfigure}{0.48\textwidth}
    \centering
    \includegraphics[width=\linewidth]{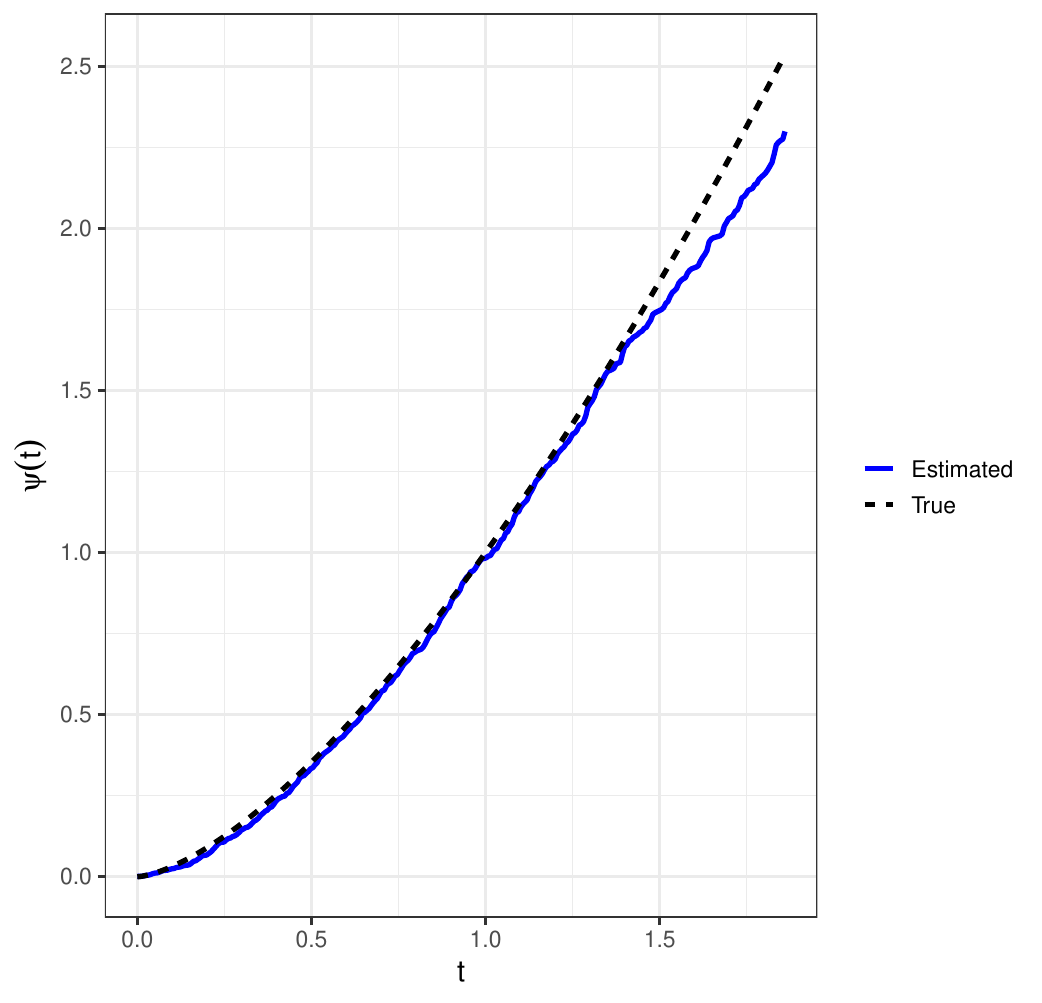}
    \caption{Results for the estimator of $\psi$ for model b).}
    \label{fig:psi-b}
\end{subfigure}

\caption{Simulation results for the identified function $\psi$.}
\label{fig:estimatorspsi}
\end{figure}


\begin{figure}[H]
\centering

\begin{subfigure}{0.48\textwidth}
    \centering
    \includegraphics[width=\linewidth]{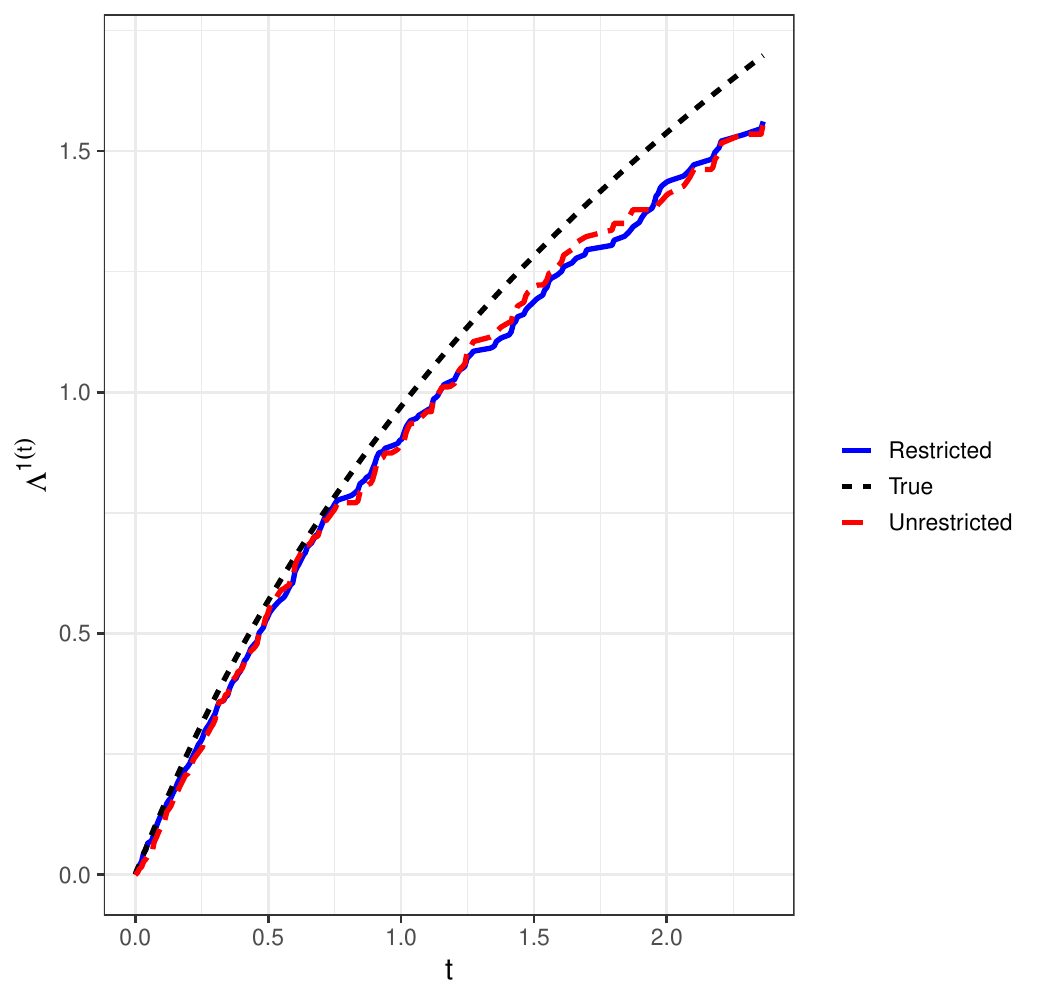}
    \caption{Results for the estimator of $\Lambda^1$ for model a).}
    \label{fig:lama}
\end{subfigure}
\hfill
\begin{subfigure}{0.48\textwidth}
    \centering
    \includegraphics[width=\linewidth]{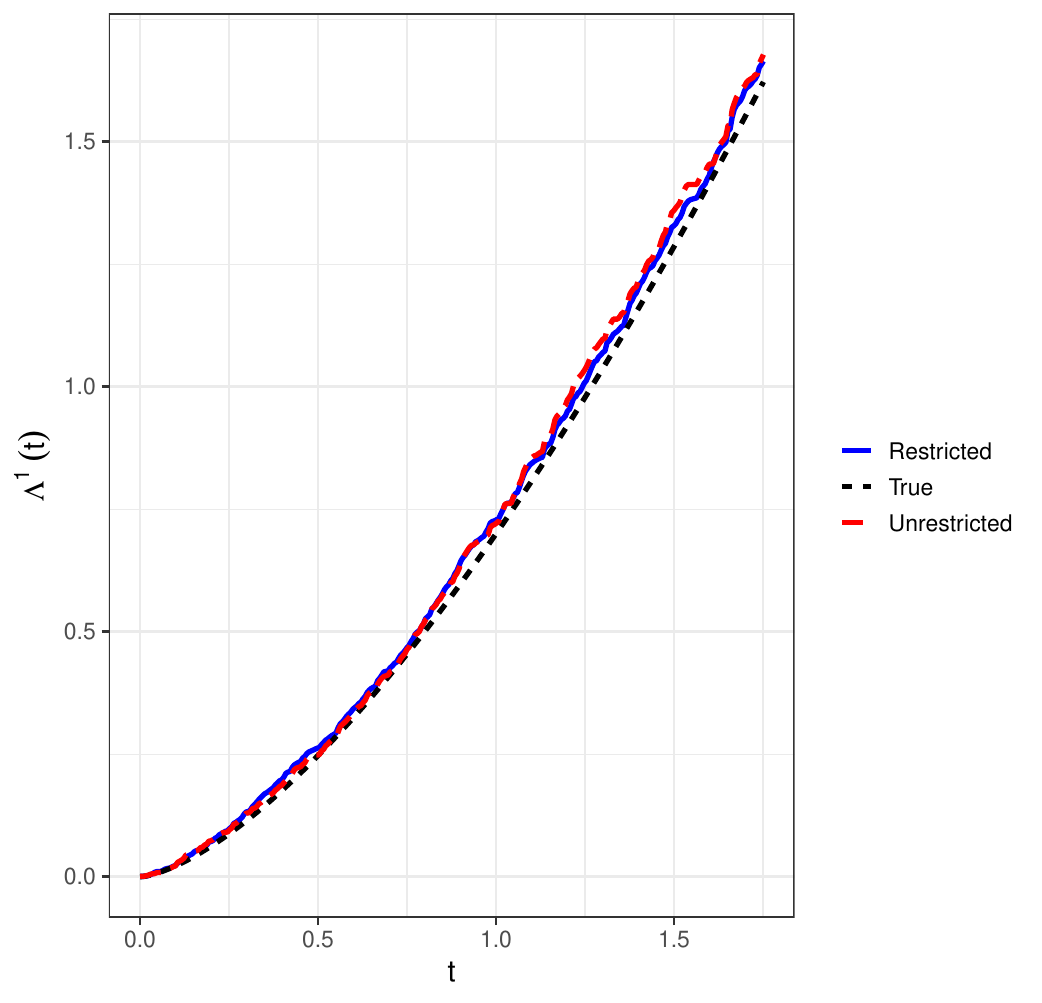}
    \caption{Results for the estimator of $\Lambda^1$ for model b).}
    \label{fig:lamb}
\end{subfigure}

\caption{Simulation results for the identified function $\Lambda^1$.}
\label{fig:estimatorslambda1}
\end{figure}


\begin{figure}[H]
\centering

\begin{subfigure}{0.48\textwidth}
    \centering
    \includegraphics[width=\linewidth]{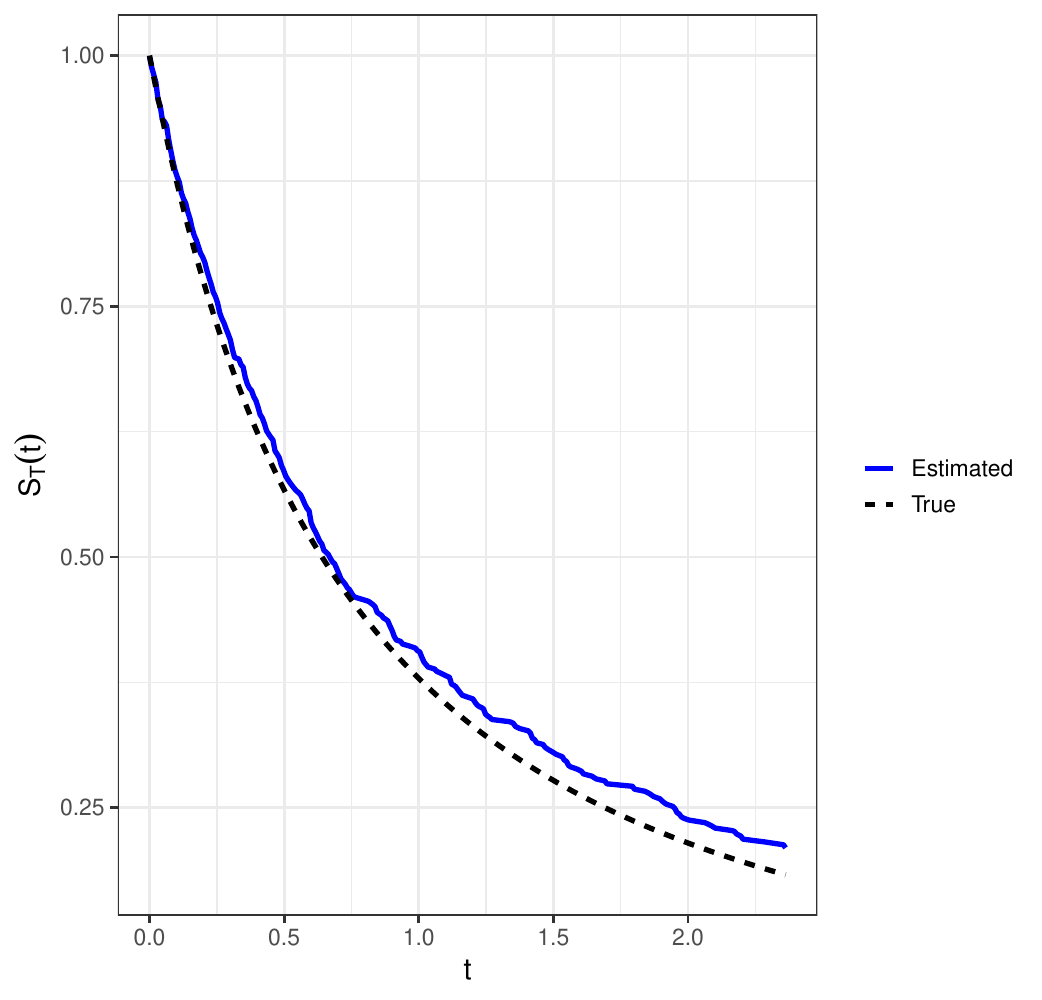}
    \caption{Results for the estimator of  $S_T$ for model a).}
    \label{fig:STa}
\end{subfigure}
\hfill
\begin{subfigure}{0.48\textwidth}
    \centering
   \includegraphics[width=\linewidth]{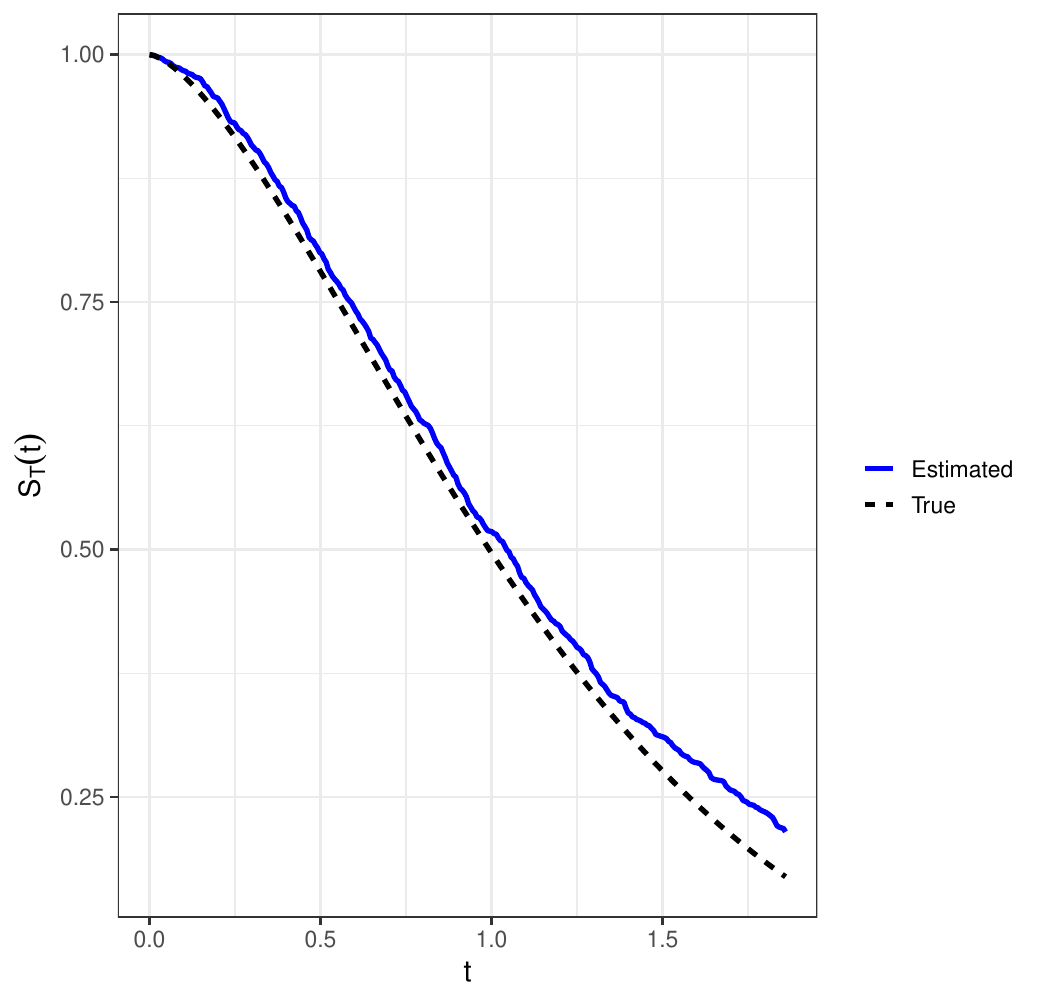}
    \caption{Results for the estimator of  $S_T$ for model b).}
    \label{fig:STb}
\end{subfigure}

\caption{Simulation results for the survival estimators of $S_T$.}
\label{fig:estimatorsST}
\end{figure}


\begin{figure}[H]
\centering

\begin{subfigure}{0.48\textwidth}
    \centering
    \includegraphics[width=\linewidth]{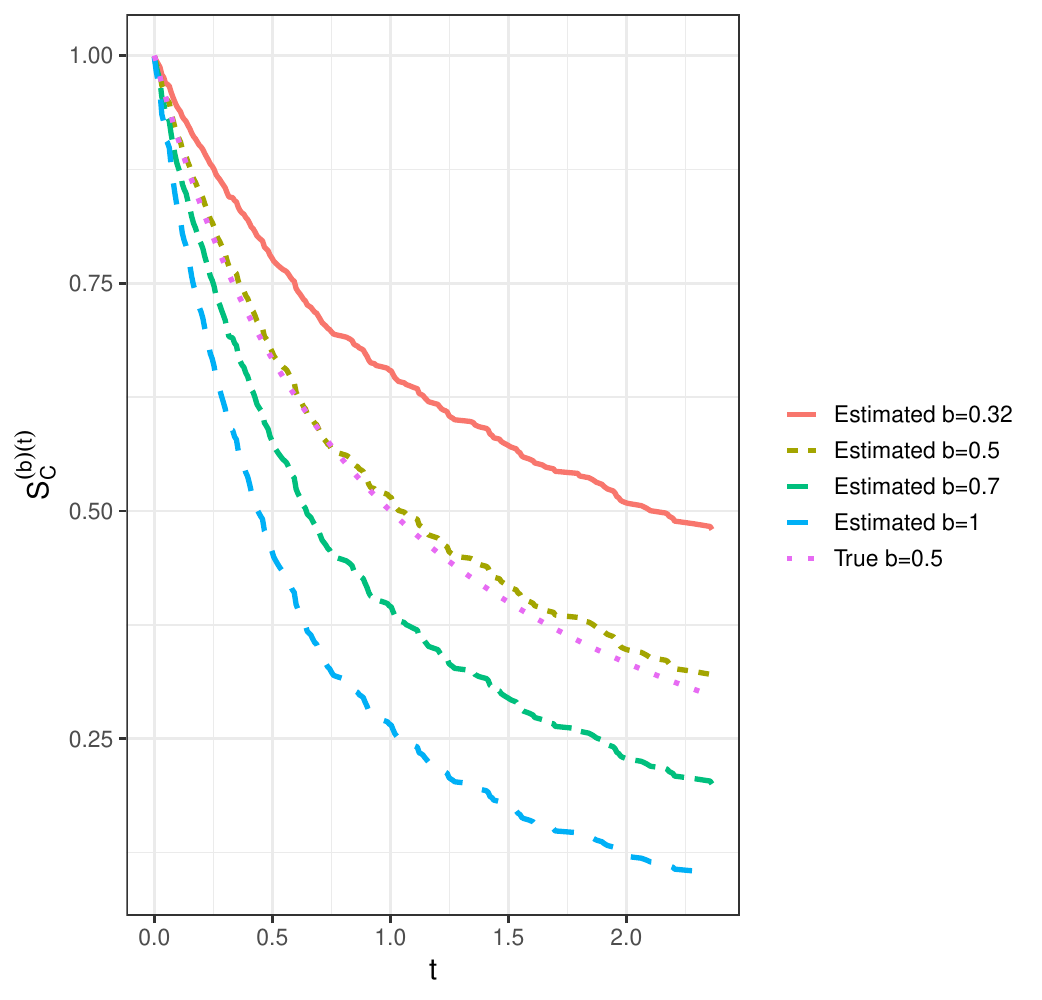}
    \caption{Results for the estimator of $S_C$ for model a).}
    \label{fig:SCa}
\end{subfigure}
\hfill
\begin{subfigure}{0.48\textwidth}
    \centering
    \includegraphics[width=\linewidth]{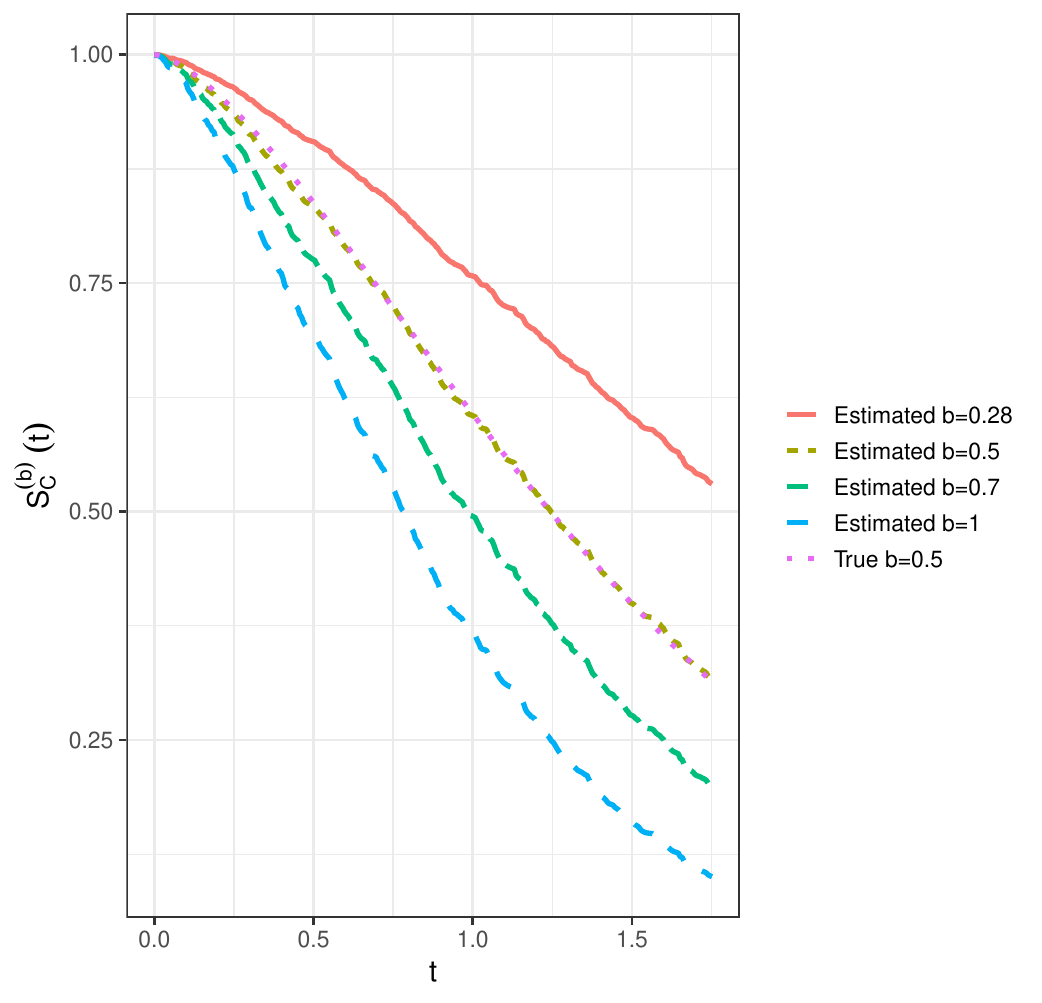}
    \caption{Results for the estimator of $S_C$ for model b).}
    \label{fig:Scb}
\end{subfigure}

\caption{Simulation results for the survival estimators of $S_C$.}
\label{fig:estimatorsSC}
\end{figure}


\begin{figure}[H]
\centering
\includegraphics[width =14cm, angle=0]{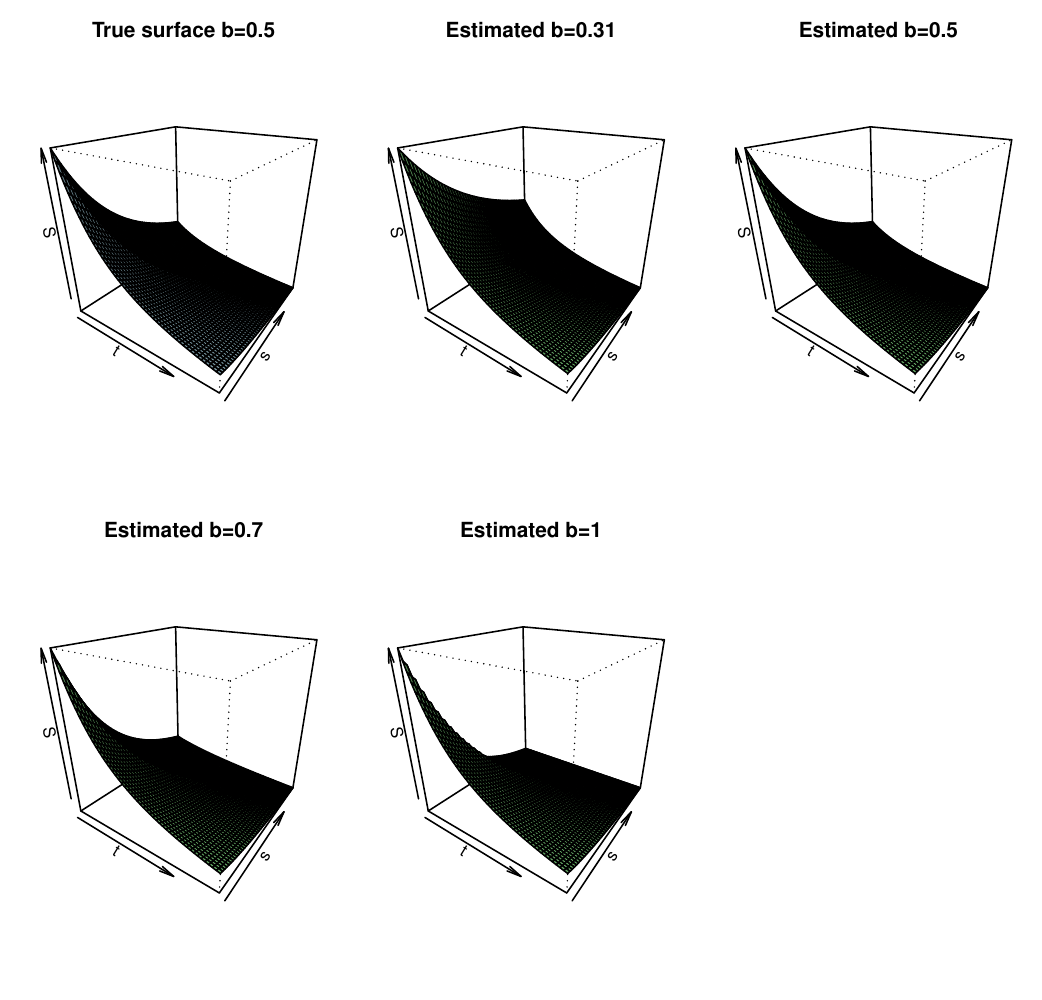}
\caption{Results for the  estimators of the joint survival function $S_{T,C}$ for model a).}\label{fig:jla}
\end{figure}


\begin{figure}[H]
\centering
\includegraphics[width =14cm, angle=0]{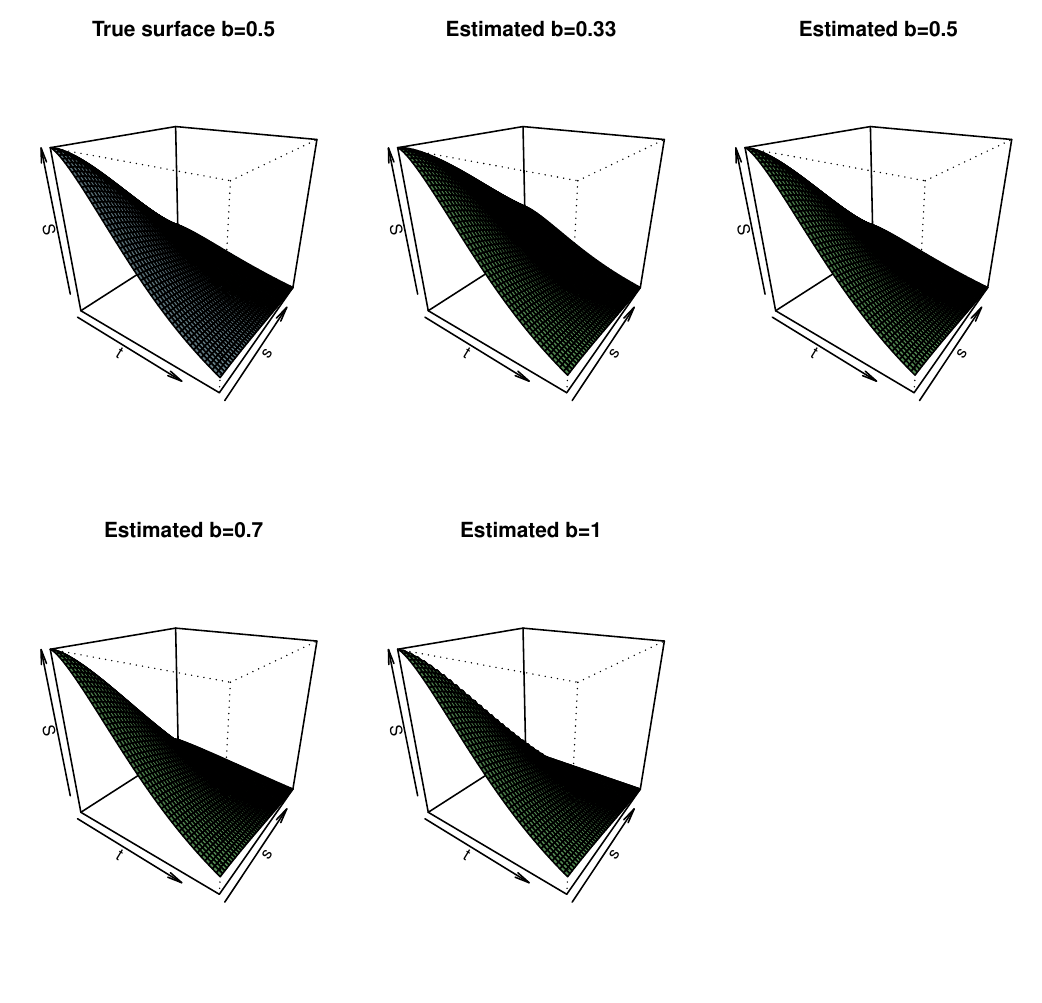}
\caption{Results for the  estimators of the joint survival function $S_{T,C}$ for model b).}\label{fig:jlb}
\end{figure}


\begin{table}[H]
\centering
\caption{Finite-sample performance of the estimator of the event loading "$a$" under the two simulation models.}
\label{tab:a-montecarlo}
\begin{tabular}{c|ccc|ccc}
\hline
& \multicolumn{3}{c|}{Model (a)} & \multicolumn{3}{c}{Model (b)}\\
\cline{2-7}
$n$ & Bias & Variance & MSE & Bias & Variance & MSE\\
\hline
50 & $0.002980$ & $0.004054$  & $0.004059$ & $0.000900$ & $0.001979$ & $0.001977$ \\
100  & $-0.001800$  & $0.0020572$ & $0.002058$ & $-0.000754$  & $0.000410$ & $0.000410$ \\
500 & $\mathbf{-0.000166}$ & $\mathbf{0.000429}$ & $\mathbf{0.000429}$ & $\mathbf{0.000437}$ & $\mathbf{0.000237}$ & $\mathbf{0.000237}$ \\
\hline
\end{tabular}
\end{table}

\begin{table}[H]
\centering
\caption{Integrated squared error (ISE) for the estimation of the identifiable quantities under the two simulation models.}
\label{tab:ISE-maxstable}
\begin{tabular}{c|cccc}
\hline
\multicolumn{5}{c}{Model (a)}\\
\hline
$n$
& $\psi$
& $\Lambda^1$ (unrestricted)
& $\Lambda^1$ (restricted)
& $S_T$\\
\hline
50 & $0.058320$ & $0.118951$ & $0.056024$ & $0.003726$
\\
100 & $0.005492$ & $0.002490$ & $ 0.005468$ & $0.000501$
\\
500 & $\mathbf{0.001973}$ & $\mathbf{0.001695}$ & $\mathbf{ 0.000724}$ & $\mathbf{0.000142}$
\\
\hline
\multicolumn{5}{c}{}\\
\hline
\multicolumn{5}{c}{Model (b)}\\
\hline
$n$
& $\psi$
& $\Lambda^1$ (unrestricted)
& $\Lambda^1$ (restricted)
& $S_T$\\
\hline
50 & $ 0.084297$  & $ 0.050385$ & $0.049474$ & $0.004932$
\\
100 & $0.008617$ & $0.003974$ & $0.002078$ & $0.000599$
\\
500 & $\mathbf{0.000827}$ & $\mathbf{0.000577}$   & $\mathbf{0.000508}$   & $\mathbf{0.000165}$
\\
\hline
\end{tabular}
\end{table}
The simulation results for the two max-stable specifications are presented
in Figures~\ref{fig:estimatorspsi}--\ref{fig:jlb}. Figure~\ref{fig:estimatorspsi} shows the estimation
of the common clock $\Psi$ for both models, while
Figure~\ref{fig:estimatorslambda1} compares the unrestricted and restricted estimators
of $\Lambda^1$. In both specifications, the estimated common clock closely follows its theoretical counterpart, and the restricted estimator of
$\Lambda^1$ remains close to the true reduction.
The unrestricted estimator exhibits a larger estimation error than the restricted estimator, as reflected by the higher ISE values reported in Table~\ref{tab:ISE-maxstable}. Moreover, the ISE decreases as the sample size increases, for both estimators, confirming the consistency of the proposed procedures.  This illustrates the stabilizing effect of exploiting the max-stable restriction $\Lambda^1=a\Psi$. Figure~\ref{fig:estimatorsST} compares the estimated event-time survival function $\widetilde S_T$ with its theoretical counterpart. The estimated survival
functions closely follow the true curves under both specifications, confirming that the good performance of the restricted procedure is preserved when the
operational clock changes from Model~(a) to Model~(b).\\
\ \\
Figure~\ref{fig:estimatorsSC} illustrates the partial identification of the censoring survival
function under the binary observation scheme. The observed data determine a family of compatible censoring
survival functions
\[
S_C^{(b)}(t)=\exp\{-b\widetilde\Psi_n(t)\},
\qquad b\in[1-\widetilde a_n,1].
\]
The figure shows that the true censoring survival function, corresponding to $b=0.5$, is contained within this estimated family. This result provides
a numerical illustration of the sharp partial identification established theoretically: the proposed estimator recovers the range of censoring
survival functions compatible with the observed binary law, rather than selecting an arbitrary point estimate for the unidentified censoring distribution.\\
\ \\
Figures~\ref{fig:jla} and~\ref{fig:jlb} display the joint survival surfaces for the two max-stable
specifications and different values of the unidentified parameter $b$. The
surfaces vary with $b$ while sharing the same observable quantities $\Psi$
and $a$, and hence the same binary law of $(Z,\delta)$. The true surface,
corresponding to $b=0.5$, is one member of the admissible family, illustrating
the partial identification of the latent joint survival distribution. \\
\ \\
Finally, Table~\ref{tab:a-montecarlo} reports a Monte Carlo study based on
$M=1000$ independent replications and sample sizes
$n\in\{50,100,500\}$. The bias, variance, and MSE of
$\widehat a_n$ decrease as the sample size increases, under both clock
specifications. This behaviour is consistent with the fact that $\widehat a_n=n^{-1}\sum_{i=1}^n\delta_i$ is the empirical mean of a Bernoulli variable with parameter $a$ and therefore satisfies the usual
$\sqrt n$-consistency. The similar results obtained for
the two operational clocks confirm that the estimation of $a$ is robust to the choice of the clock function.
\subsection{Real data  application}\label{sub2}
To illustrate the practical relevance of the proposed methodology, we consider
the \texttt{veteran} dataset from the \texttt{survival} package in \textsf{R}.
This dataset comes from a randomized clinical trial on patients with advanced
lung cancer. The study records the survival time of patients together with the occurrence of the event of interest and possible right censoring.
The observed data naturally fit the binary observation framework considered in this paper. For each patient, we define \[Z_i=T_i\wedge C_i, \qquad \delta_i=\mathbf{1}_{\{T_i\le C_i\}},\]
where \(T_i\) denotes the latent event time, \(C_i\) the censoring time, 
\(Z_i\) the observed follow-up time, and \(\delta_i\) the event indicator.
After removing incomplete observations, the dataset contains \(n=137\)
patients, with an empirical event probability
\[
\widetilde a_n=\frac{1}{n}\sum_{i=1}^n\delta_i
=0.9343.
\]
Before applying the restricted max-stable estimation procedure, we first assess
whether the observable binary law is compatible with the max-stable
common-clock restriction. According to Theorem~\ref{thm:compatibility-main},
this restriction is equivalent to the independence condition $Z\perp\delta$, 
or, equivalently, to the constancy of the observable cause ratio
\[
q(t)=\mathbb P(\delta=1\mid Z=t).
\]
We therefore apply the permutation-based Kolmogorov--Smirnov test introduced in
Section~\ref{sec:testing}. The test statistic is based on
\[
\mathbb T_n(t)
=
\sqrt n
\left\{
\widetilde F_{1,n}(t)
-
\widetilde a_n\widetilde F_n(t)
\right\}.
\]
For the veteran data, the resulting test statistic is
\[
D_n=1.264,
\]
with permutation p-value
\[
p_{\mathrm{perm}}=0.065.
\]
Hence, at the 5\% significance level, the observable max-stable compatibility
hypothesis is not rejected. This supports the use of the restricted estimator
based on the proportionality relation
\[
\Lambda^1(t)=a\Lambda^{12}(t).
\]
We then estimate the common clock using the structural Nelson--Aalen estimator
\(\widetilde\Psi(t)=\widetilde\Lambda^{12}(t)\), and obtain the restricted estimator of the event reduction
\(\widetilde\Lambda^{1,R}(t)=\widetilde a_n\widetilde\Psi(t)\).
Since the censoring loading is not point identified under the max-stable model,
we consider the identified set
\[
b\in[1-\widetilde a_n,1]
=
[0.0657,1],
\]
and construct the associated family of compatible censoring survival functions
\[
\widetilde S_C^{(b)}(t)
=
\exp\{-b\widetilde\Psi(t)\},
\qquad b\in[1-\widetilde a_n,1].
\]
Finally, we estimate the event survival function
\[
\widetilde S_T(t)
=
\exp\{-\widetilde a_n\widetilde\Psi(t)\},
\]
as well as the family of compatible joint survival functions
\[
\widetilde S_{T,C}^{(b)}(t,s)
=
\exp\left(
-\widetilde a_n\widetilde\Psi(t)
-b\widetilde\Psi(s)
+(\widetilde a_n+b-1)
\widetilde\Psi(t\wedge s)
\right).
\]

\begin{figure}[H]
\centering

\begin{subfigure}{0.48\textwidth}
    \centering
    \includegraphics[width=\linewidth]{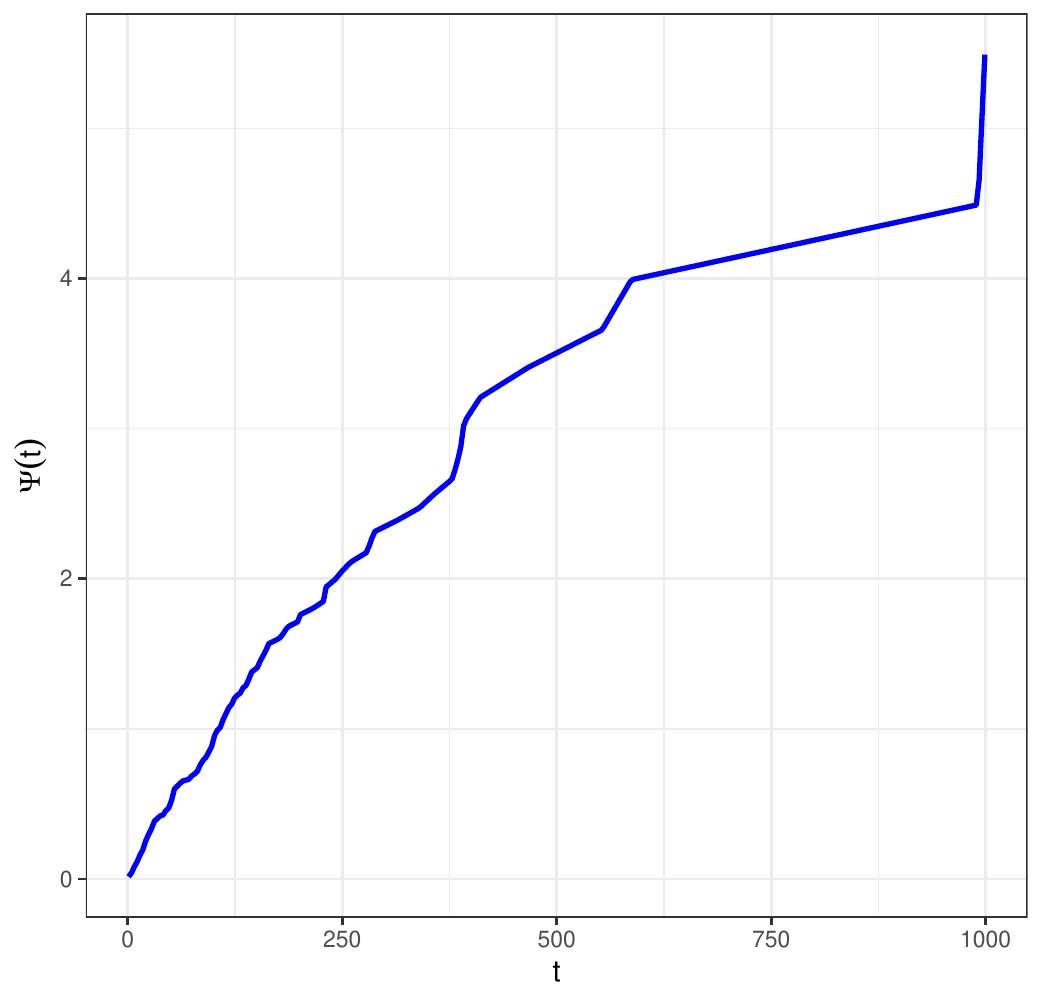}
    \caption{Results for the estimator of $\psi$.}
    \label{fig:psiA}
\end{subfigure}
\hfill
\begin{subfigure}{0.48\textwidth}
    \centering
    \includegraphics[width=\linewidth]{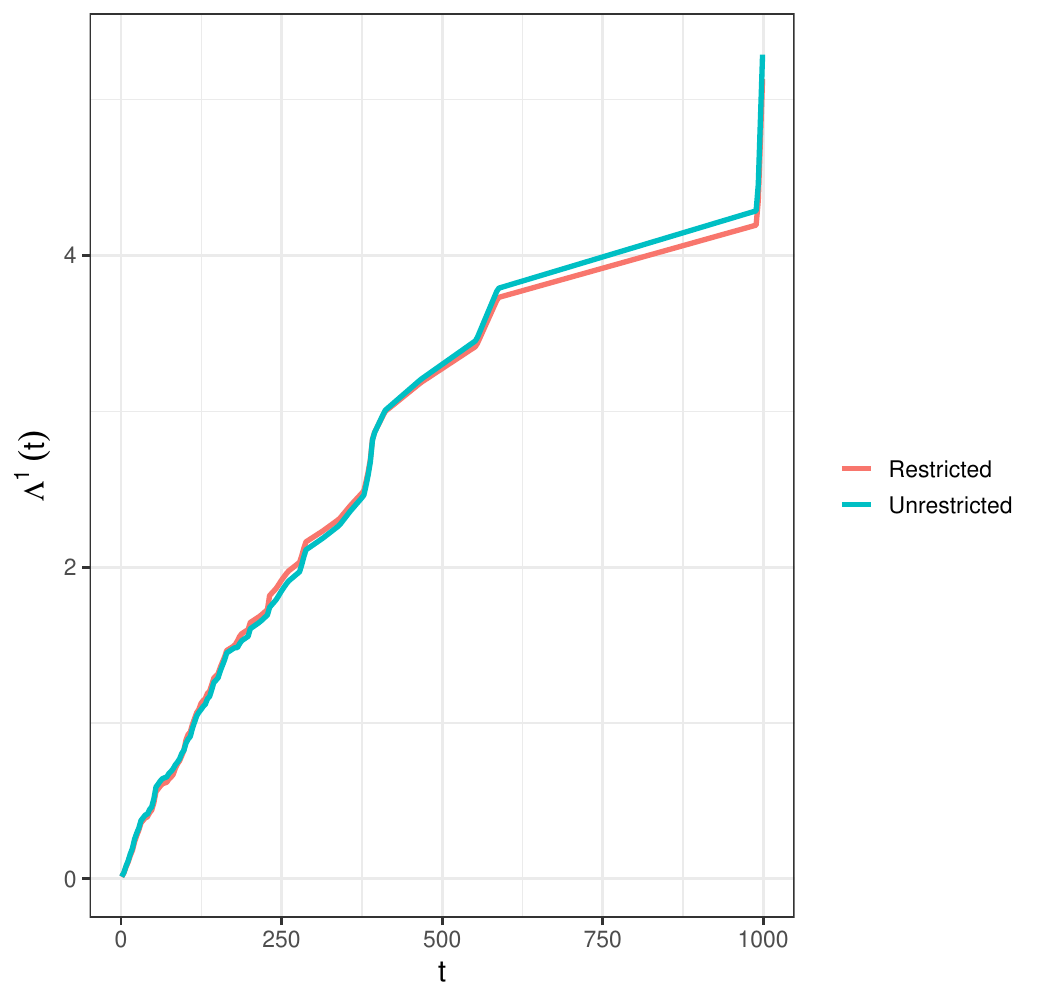}
    \caption{Results for the estimator of $\Lambda^1$.}
    \label{fig:lambda2}
\end{subfigure}

\caption{Simulation results for the identified quantities for real data.}
\label{fig:estimatorsdatapslambda1}
\end{figure}

\begin{figure}[H]
\centering

\begin{subfigure}{0.48\textwidth}
    \centering
    \includegraphics[width=\linewidth]{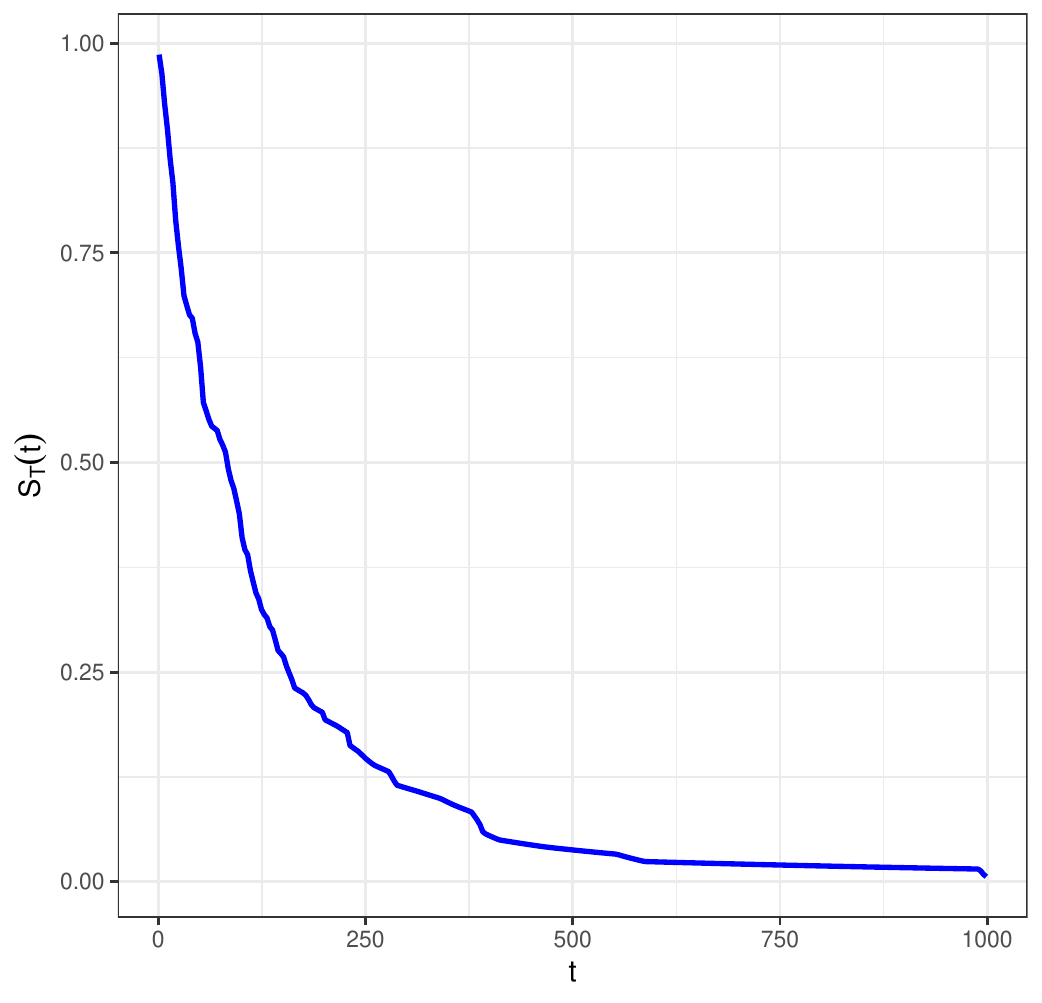}
    \caption{Results for the estimator of $S_T$.}
    \label{fig:psiB}
\end{subfigure}
\hfill
\begin{subfigure}{0.48\textwidth}
    \centering
    \includegraphics[width=\linewidth]{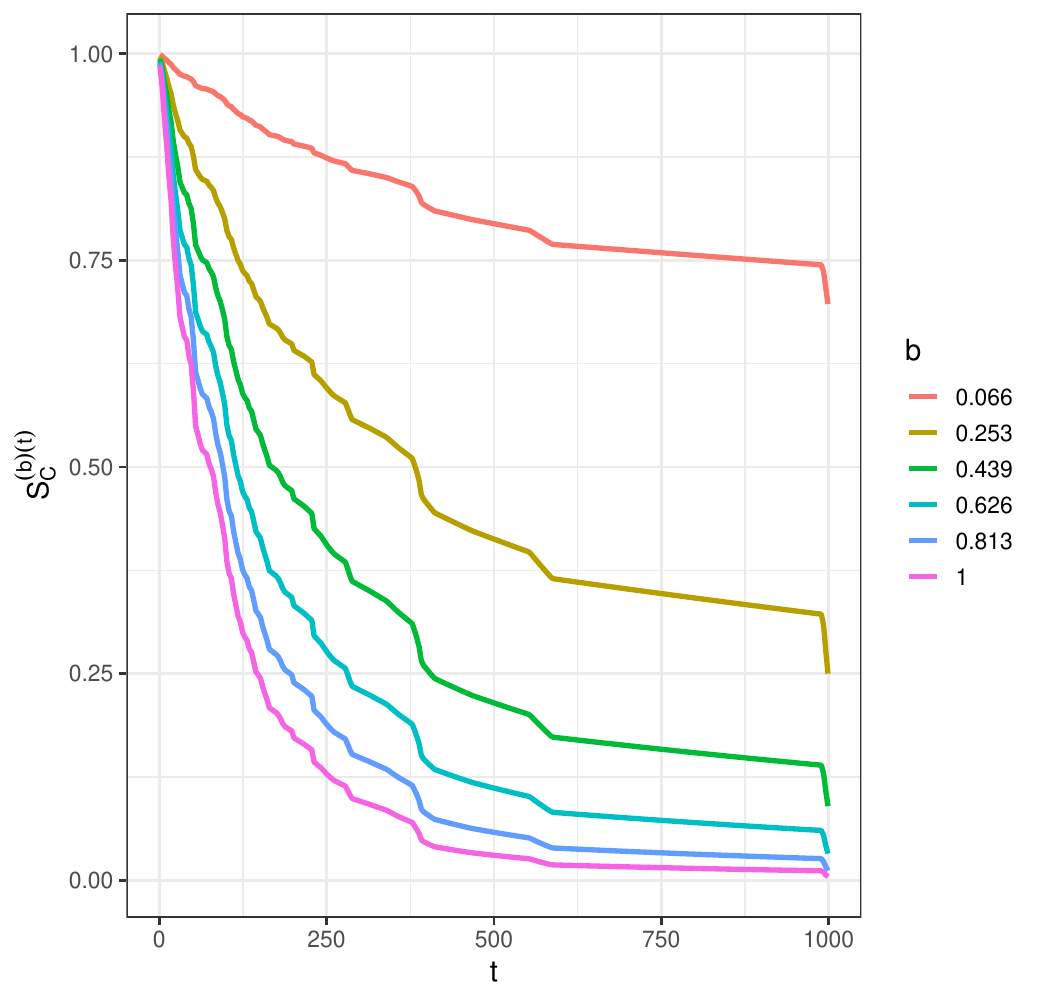}
    \caption{Results for the estimator of $S_C$.}
    \label{fig:lambda1}
\end{subfigure}

\caption{Simulation results for the survival estimators for real data.}
\label{fig:estimatorsdataSTSC}
\end{figure}

\begin{figure}[H]
\centering
\includegraphics[width =14cm, angle=0]{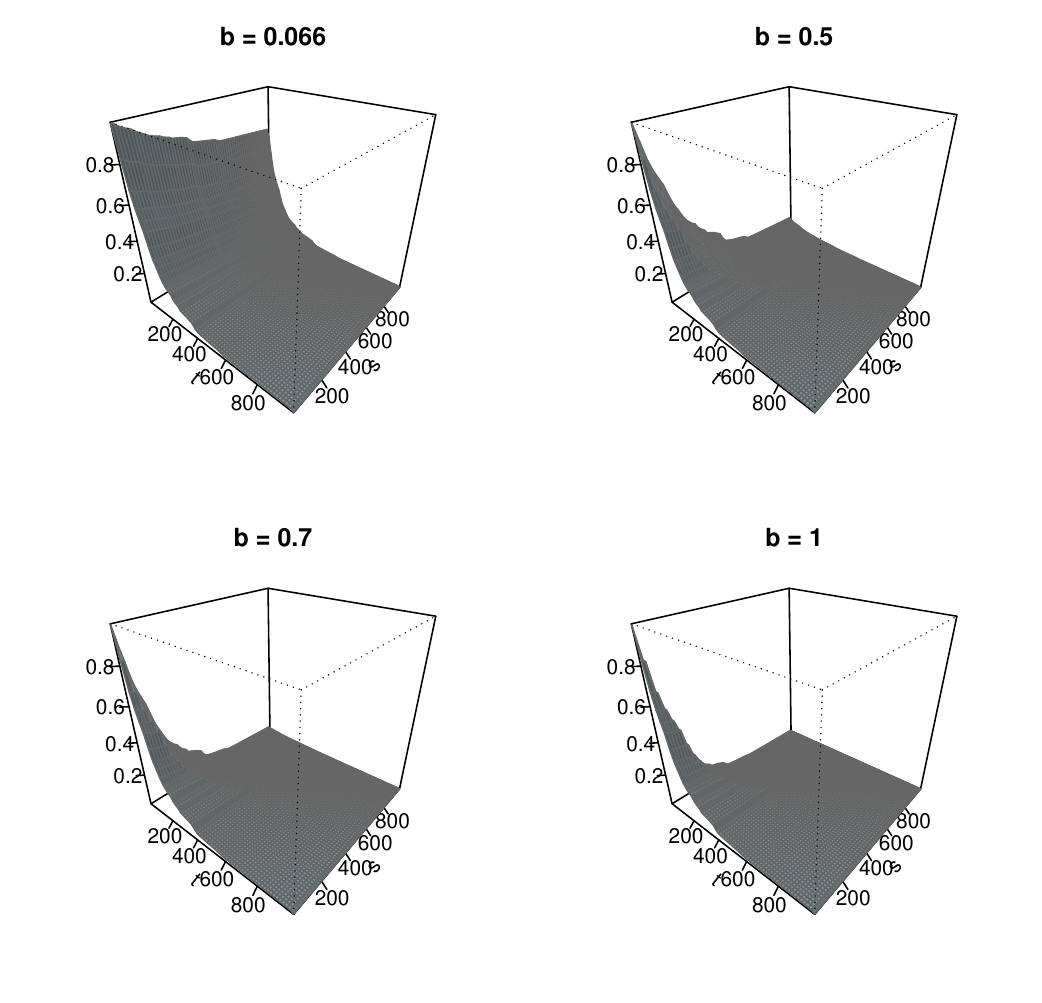}
\caption{Results for the  estimators of the joint survival function $S_{T,C}$ for real data.}\label{fig:jldata}
\end{figure}

The estimated common clock $\widetilde\Psi$ and the restricted estimator of $\Lambda^1$ provide a stable description of the observable event mechanism, with the corresponding estimate of the event survival function $\widetilde S_T$ displaying the expected decreasing behaviour.   The estimated censoring survival functions illustrate the remaining
non-identification of the censoring mechanism. Since
$b\in[0.0657,1]$, the data are compatible with a broad family of censoring survival functions rather than with a unique censoring distribution. The
corresponding joint survival surfaces similarly vary with $b$, showing that the latent dependence structure between $T$ and $C$ cannot be uniquely
recovered from the binary observations.
Overall, the real-data analysis illustrates the main practical implication of the proposed framework: the event-time distribution and the common clock
are point identified under observable compatibility, whereas the censoring distribution and the joint survival distribution remain only partially
identified.

\section{Conclusion and Discussion}
\label{sec:conclusion}

The results complement copula-graphic approaches \citep{ZhengKlein1995,RivestWells2001}, parametric copula methods for dependent censoring \citep{CzadoVanKeilegom2023}, and generalized Marshall--Olkin estimation procedures \citep{EscobarBachHelali2024,Helali2025}. The aim here is different: to give the exact information content of binary and complete records within the max-stable reduction class.
Max-stability strongly reduces the dimension of a dependent-censoring model, but it cannot create information absent from the observation scheme. In the complete three-category experiment, the simultaneous category identifies the common-shock component. In the binary experiment, this component is merged with strict failures; the correct output is a point estimate for $(\Psiop,a)$ and sharp identified sets for the latent quantities depending on $b$. The efficient max-stable restricted estimator is efficient for the point-identified observable target, not for a latent completion that the binary law cannot identify. This distinction is the central methodological message of the paper.

\appendix

\section{Proofs}
\label{sec:proofs}

\subsection{Proof of Proposition~\ref{prop:submeasures} }
We first derive the law of \(Z=T\wedge C\). Since \(\Pp(Z>t)=\Pp(T>t,C>t)=H(t)=e^{-\Lambda^{12}(t)}\), we obtain \(f_Z(t)=-H'(t)=H(t)(\Lambda^{12})'(t)\), which proves \eqref{eq:fz}.
To identify the different causes of failure, consider first \(t<s\). From \eqref{eq:joint-survival},
\[
S_{T,C}(t,s)=\exp\{-\Lambda^2(s)-\Lambda^{12}(t)+\Lambda^2(t)\},
\]
and hence \(-\partial_tS(t,s)=S(t,s)\{(\Lambda^{12})'(t)-(\Lambda^2)'(t)\}\). Letting \(s\downarrow t\) gives \(f_{Z,\Delta=1}(t)=H(t)\{(\Lambda^{12})'(t)-(\Lambda^2)'(t)\}\). Since \(\Gamma'=(\Lambda^1)'+(\Lambda^2)'-(\Lambda^{12})'\), this is equivalently \(f_{Z,\Delta=1}(t)=H(t)\{(\Lambda^1)'(t)-\Gamma'(t)\}\), proving \eqref{eq:f1}.
Similarly, for \(s<t\), we have
\[
S_{T,C}(t,s)=\exp\{-\Lambda^1(t)-\Lambda^{12}(s)+\Lambda^1(s)\},
\]
so that \(-\partial_sS(t,s)=S(t,s)\{(\Lambda^{12})'(s)-(\Lambda^1)'(s)\}\). Letting \(s\uparrow t\) yields \(f_{Z,\Delta=0}(t)=H(t)\{(\Lambda^{12})'(t)-(\Lambda^1)'(t)\}\). Equivalently, using the definition of \(\Gamma\), \(f_{Z,\Delta=0}(t)=H(t)\{(\Lambda^2)'(t)-\Gamma'(t)\}\), proving \eqref{eq:f0}.
Finally, the three complete-category submeasures partition the law of \(Z\), so \(f_{Z,\Delta=2}=f_Z-f_{Z,\Delta=1}-f_{Z,\Delta=0}\). Using \eqref{eq:fz}--\eqref{eq:f0},
\[
\begin{aligned}
f_{Z,\Delta=2}(t)
&=H(t)\Bigl[(\Lambda^{12})'
-\{(\Lambda^{12})'-(\Lambda^2)'\}
-\{(\Lambda^{12})'-(\Lambda^1)'\}\Bigr]\\
&=H(t)\{(\Lambda^1)'+(\Lambda^2)'-(\Lambda^{12})'\}
=H(t)\Gamma'(t),
\end{aligned}
\]
which proves \eqref{eq:f2}.
Finally, since \(\{\delta=1\}=\{\Delta=1\}\cup\{\Delta=2\}\) and \(\{\delta=0\}=\{\Delta=0\}\), with disjoint union, we obtain \(f_{Z,\delta=1}(t)=f_{Z,\Delta=1}(t)+f_{Z,\Delta=2}(t)=H(t)(\Lambda^1)'(t)\) and \(f_{Z,\delta=0}(t)=f_{Z,\Delta=0}(t)=H(t)\{(\Lambda^{12})'(t)-(\Lambda^1)'(t)\}\), which is exactly \eqref{eq:binary-submeasures}.
\subsection{Proof of Theorem~\ref{thm:maxstab}  }
{The argument is the bivariate event--censoring specialization of the common-clock characterization for max-stable Sibuya survival copulas in \citet{GueyeQuintos2025}.
We first prove sufficiency. Suppose that
\(\Lambda^1=\lambda_1\psi\), \(\Lambda^2=\lambda_2\psi\), and
\(\Lambda^{12}=\lambda_{12}\psi\), for a continuous, strictly increasing,
unbounded clock \(\psi\) and positive constants
\(\lambda_1,\lambda_2,\lambda_{12}\). Then
\[
\Gamma=(\lambda_1+\lambda_2-\lambda_{12})\psi.
\]
For \(u=e^{-\lambda_1\psi(t)}\) and \(v=e^{-\lambda_2\psi(s)}\), we have
\(t_u=\psi^{-1}(-\log u/\lambda_1)\) and
\(s_v=\psi^{-1}(-\log v/\lambda_2)\). Substituting these expressions into
\eqref{eq:surv-copula} gives
\[
\Cop(u,v)
=
uv\exp\!\left\{
(\lambda_1+\lambda_2-\lambda_{12})
\min\!\left(
\frac{-\log u}{\lambda_1},
\frac{-\log v}{\lambda_2}
\right)
\right\}.
\]
Replacing \((u,v)\) by \((u^r,v^r)\) simply multiplies both arguments of the
minimum by \(r\). Hence
\(\log\Cop(u^r,v^r)=r\log\Cop(u,v),
\)
which is equivalent to
\(\Cop(u^r,v^r)=\Cop(u,v)^r.\)
Therefore the copula is max-stable.\\
\ \\
We now prove the converse. Assume that \(\Cop\) is max-stable and that
\(\Gamma\not\equiv0\).
Set
$$
R_1(x)=(\Lambda^1)^{-1}(x),
\qquad
R_2(y)=(\Lambda^2)^{-1}(y),
\qquad x,y>0.
$$
In logarithmic marginal coordinates $u=e^{-x}$ and $v=e^{-y}$, the copula satisfies
$$
-\log \Cop(e^{-x},e^{-y})
=
x+y-\Gamma\{R_1(x)\wedge R_2(y)\}.
$$
Max-stability is therefore equivalent to the homogeneity identity
\[
\Gamma\{R_1(rx)\wedge R_2(ry)\}
=
r\,\Gamma\{R_1(x)\wedge R_2(y)\},
\qquad x,y>0,\ r>0. \tag{MS}
\]
Fix $x>0$ and $r>0$. Since $R_2(y)\to\infty$ as $y\to\infty$, one may choose $y$ large enough, for this fixed $r$, such that
$$
R_2(y)>R_1(x),
\qquad
R_2(ry)>R_1(rx).
$$
Then (MS) gives $$\Gamma(R_1(rx))=r\Gamma(R_1(x)).
$$ Hence $h_1(x)=\Gamma(R_1(x))$ is homogeneous of degree one: $h_1(rx)=rh_1(x)$. Thus $h_1(x)=c_1x$ and $\Gamma(t)=c_1\Lambda^1(t)$. By interchanging the roles of $1$ and $2$, we obtain $\Gamma(t)=c_2\Lambda^2(t)$. Therefore $\Lambda^1$ and $\Lambda^2$ are proportional, and then $\Lambda^{12}=\Lambda^1+\Lambda^2-\Gamma$ is proportional to the same clock.


\subsection{Proof of Proposition \ref{prop:copula-pickands}  }
Under the normalization \eqref{eq:normalized-model}, we have \(
\Gamma(t)=(a+b-1)\Psi(t)=\theta\Psi(t),\)
where \(\theta=a+b-1\ge0\).
Moreover, $\Psi(t)=\Lambda^{12}(t).$
For the marginal reductions,
\[
u=e^{-\Lambda^1(t)}
=e^{-a\Psi(t)},
\qquad
v=e^{-\Lambda^2(s)}
=e^{-b\Psi(s)}.
\]
Hence,
\(
\Psi(t)=-\frac{1}{a}\log u,
\Psi(s)=-\frac{1}{b}\log v .
\)
Using the survival copula representation
\[
\Cop(u,v)
=
uv\exp\{\Gamma(t_u\wedge s_v)\},
\]
we obtain
\(
\Cop(u,v)
=
uv
\exp\left\{
\theta\Psi(t_u\wedge s_v)
\right\}.
\)
Since \(\Psi\) is increasing,
\(
\Psi(t_u\wedge s_v)
=
\min\{\Psi(t_u),\Psi(s_v)\}.
\)
Therefore,
\(
\Psi(t_u\wedge s_v)
=
\min\left\{-\frac{\log u}{a},
-\frac{\log v}{b}\right\}.
\)
Consequently,
\[
\Cop_{a,b}(u,v)
=
uv
\exp\left\{
\theta
\min\left(
-\frac{\log u}{a},
-\frac{\log v}{b}
\right)
\right\}.
\]
Using
$e^{c\min(x,y)}
=
\min(e^{cx},e^{cy}),$
with
$c=\theta,$
we get
\(
\Cop_{a,b}(u,v)
=
uv
\min\left\{
u^{-\theta/a},
v^{-\theta/b}
\right\},
\)
which proves \eqref{eq:MO-copula}.
We now derive the stable-tail dependence function. By definition,
\[
\ell_{a,b}(x,y)
=
-\log\Cop_{a,b}(e^{-x},e^{-y}).
\]
Substituting \((u,v)=(e^{-x},e^{-y})\) into
\eqref{eq:MO-copula}, we obtain
\(
\Cop_{a,b}(e^{-x},e^{-y})
=
e^{-x-y}
\min\left\{
e^{\theta x/a},
e^{\theta y/b}
\right\}.
\)
Taking minus the logarithm gives
\(
\ell_{a,b}(x,y)
=
x+y
-
\min\left\{
\frac{\theta x}{a},
\frac{\theta y}{b}
\right\},
\)
which proves \eqref{eq:ell-ab}.
Finally, the Pickands dependence function is defined by
\[
A_{a,b}(w)=\ell_{a,b}(w,1-w),
\qquad w\in[0,1].
\]
Hence,
\(
A_{a,b}(w)
=
1-
\min\left\{
\frac{\theta w}{a},
\frac{\theta(1-w)}{b}
\right\}.
\)
The breakpoint is obtained when the two affine functions inside the minimum
are equal:
\[
\frac{\theta w}{a}
=
\frac{\theta(1-w)}{b}.
\]
Since \(\theta>0\) in the genuinely dependent case, this is equivalent to
\(
bw=a(1-w),
\)
and therefore
\(
w^*(a,b)=\frac{a}{a+b}.
\)
Thus \(A_{a,b}\) is the minimum of two affine functions and is consequently
piecewise affine, with breakpoint \(w^*(a,b)\).

\subsection{Proof of Theorem 4.1 } 
The equality of the binary laws implies
\[
\widetilde{f}_{Z,\delta=1}(t)=f_{Z,\delta=1}(t)
\quad\text{and}\quad
\widetilde{f}_{Z,\delta=0}(t)=f_{Z,\delta=0}(t).
\]
Since
\(
 f_{Z,\delta=1}(t)=H(t)(\Lambda^1)'(t),
 f_{Z,\delta=0}(t)=H(t)\bigl\{(\Lambda^{12})'(t)-(\Lambda^1)'(t)\bigr\},
\)
with \(H(t)=e^{-\Lambda^{12}(t)}\), we obtain
\[
e^{-\widetilde{\Lambda}^{12}(t)}
\bigl(\widetilde{\Lambda}^{1}\bigr)'(t)
=
e^{-\Lambda^{12}(t)}
\bigl(\Lambda^{1}\bigr)'(t),
\]
and
\[
e^{-\widetilde{\Lambda}^{12}(t)}
\left\{
\bigl(\widetilde{\Lambda}^{12}\bigr)'(t)
-
\bigl(\widetilde{\Lambda}^{1}\bigr)'(t)
\right\}
=
e^{-\Lambda^{12}(t)}
\left\{
(\Lambda^{12})'(t)
-
(\Lambda^{1})'(t)
\right\}.
\] 
Adding the two equalities yields
\(
e^{-\widetilde{\Lambda}^{12}(t)}
\bigl(\widetilde{\Lambda}^{12}\bigr)'(t)
=
e^{-\Lambda^{12}(t)}
\bigl(\Lambda^{12}\bigr)'(t).
\)
Since
\[
\frac{d}{dt}\!\left(-e^{-\Lambda^{12}(t)}\right)
=
e^{-\Lambda^{12}(t)}(\Lambda^{12})'(t),
\]
we conclude that
\(
\frac{d}{dt}\!\left(e^{-\widetilde{\Lambda}^{12}(t)}
-e^{-\Lambda^{12}(t)}\right)=0.
\)
Hence,
\(
e^{-\widetilde{\Lambda}^{12}(t)}
-
e^{-\Lambda^{12}(t)}
=C
\)
for some constant \(C\). Since \(\widetilde{\Lambda}^{12}(0)=\Lambda^{12}(0)=0\), it follows that \(C=0\). Therefore,
\(
\widetilde{\Lambda}^{12}(t)=\Lambda^{12}(t).
\)
Substituting this identity into the first equality gives
$(\widetilde{\Lambda}^{1})'(t)
=(\Lambda^{1})'(t)$.
Using again \(\widetilde{\Lambda}^{1}(0)=\Lambda^{1}(0)=0\), we conclude that $\widetilde{\Lambda}^{1}(t)=\Lambda^{1}(t)$,
which proves the identifiability of \((\Lambda^{1},\Lambda^{12})\).\\
\ \\
Conversely, let $(\Lambda^1,\Lambda^2,\Lambda^{12})$ and 
$(\widetilde\Lambda^1,\widetilde\Lambda^2,\widetilde\Lambda^{12})$
be two admissible triplets such that $$\Lambda^{12}=\widetilde\Lambda^{12},\qquad
\Lambda^1=\widetilde\Lambda^1.$$ Then
\(
H(t)=e^{-\Lambda^{12}(t)}
=
e^{-\widetilde\Lambda^{12}(t)}
=\widetilde H(t),\)
so the two triplets induce the same distribution of the observed minimum
$Z=T\wedge C$. Moreover,
$$
\mathbb P(Z\in\dd t,\delta=1)
=
H(t)\dd\Lambda^1(t),
$$
and
$$
\mathbb P(Z\in\dd t,\delta=0)
=
H(t)\{\dd\Lambda^{12}(t)-\dd\Lambda^1(t)\}.
$$
Hence both subdistribution measures depend only on
$(\Lambda^1,\Lambda^{12})$, and therefore coincide for the two triplets.
Consequently,
\(
\mathcal L(Z,\delta)
=
\widetilde{\mathcal L}(Z,\delta).\)
This proves
$$
\mathcal L(Z,\delta)
=
\widetilde{\mathcal L}(Z,\delta)
\quad\Longleftrightarrow\quad
\Lambda^{12}=\widetilde\Lambda^{12},
\qquad
\Lambda^1=\widetilde\Lambda^1.
$$
It remains to characterize the sharp identified set for $\Lambda^2$.
Fix the identified pair $(\Lambda^1,\Lambda^{12})$, and let $\Lambda^2=L.$ The corresponding interaction reduction is
\(
\Gamma_L=\Lambda^1+L-\Lambda^{12}.
\)
Under the absolute-continuity assumptions, admissibility requires
$$
\Gamma_L'\ge0,
\qquad
(\Lambda^1)'-\Gamma_L'\ge0,
\qquad
L'-\Gamma_L'\ge0.
$$
Using
\(
\Gamma_L'
=
(\Lambda^1)'+L'-(\Lambda^{12})',
 \)
these conditions are equivalent to
$$
(\Lambda^{12})'-(\Lambda^1)'
\le L'
\le(\Lambda^{12})'
\qquad\text{a.e.},
$$
the remaining inequality
\(
(\Lambda^{12})'-(\Lambda^1)'\ge0
\)
being already implied by the nonnegativity of $\mathbb P(Z\in\dd t,\delta=0).$
Therefore the sharp identified set is exactly
$$
\mathcal I_2(\Lambda^1,\Lambda^{12})
=
\left\{
L:
L(0)=0,\;
(\Lambda^{12})'-(\Lambda^1)'
\le L'
\le(\Lambda^{12})'
\ \text{a.e.}
\right\}.
$$
Every $L$ in this set defines an admissible triplet and, by the previous equivalence, generates the same binary law. Hence the set is sharp.

\subsection{Proof of Theorem \ref{thm:binary-MS} }
Under the normalized max-stable representation
\eqref{eq:normalized-model}, we have
\(
\Lambda^1=a\Psiop,
\Lambda^{12}=\Psiop.
\)
Using the binary submeasures in
\eqref{eq:binary-submeasures}, we obtain
\[
f_{Z,\delta=1}(t)
=
H(t)(\Lambda^1)'(t),
\quad
f_{Z,\delta=0}(t)
=
H(t)\bigl\{(\Lambda^{12})'(t)-(\Lambda^1)'(t)\bigr\}.
\]
Since
\(
H(t)=e^{-\Lambda^{12}(t)}
=e^{-\Psiop(t)},
\)
\(
(\Lambda^1)'(t)=a\Psiop'(t)\) and 
\( (\Lambda^{12})'(t)=\Psiop'(t),\)
we get
\[
f_{Z,\delta=1}(t)
=
a e^{-\Psiop(t)}\Psiop'(t),
\quad
f_{Z,\delta=0}(t)
=
(1-a)e^{-\Psiop(t)}\Psiop'(t).
\]
Equivalently,
\(
\Pp(Z\in\dd t,\delta=1)
=
a e^{-\Psiop(t)}\dd\Psiop(t),
\)
and
\(
\Pp(Z\in\dd t,\delta=0)
=
(1-a)e^{-\Psiop(t)}\dd\Psiop(t),
\)
which proves \eqref{eq:binary-factorization}.
Adding the two submeasures gives the marginal law of \(Z\):
\[
\Pp(Z\in\dd t)
=
e^{-\Psiop(t)}\dd\Psiop(t).
\]
Therefore,
\(
\Pp(Z>t)=e^{-\Psiop(t)},
\)
and hence
\(
\Psiop(t)
=
-\log\Pp(Z>t).
\)
Thus the distribution of \(Z\) identifies the common clock
\(\Psiop\).
Moreover,
\[
\Pp(\delta=1)
=
\int_0^\infty
a e^{-\Psiop(t)}\dd\Psiop(t).
\]
Since
\(
\int_0^\infty e^{-\Psiop(t)}\dd\Psiop(t)
=
1,
\)
we obtain
\(
\Pp(\delta=1)=a.
\)
Hence \(a\) is identified by the binary observation.
We now show that \(b\) is not identified.
The binary observable law in
\eqref{eq:binary-factorization} depends only on
\((\Psiop,a)\) and does not involve \(b\).
Therefore, two different values of \(b\) leading to the same admissible
normalized model generate exactly the same distribution of
\((Z,\delta)\).
It remains to characterize the admissible values of \(b\).
From the constraints
\[
0<a\le1,
0<b\le1,
a+b\ge1,
\]
we obtain
\(
1-a\le b\le1.
\)
Conversely, for every
\(
b\in[1-a,1],
\)
the conditions above are satisfied. Indeed,
\(
b>0, b\le1,
\)
and
\(
a+b\ge a+(1-a)=1.
\)
Therefore every value
\[
b\in[1-a,1]
\]
defines an admissible max-stable model with the same identified pair
\((\Psiop,a)\) and consequently the same binary observable law.
Hence the sharp identified set is
\(
\mathcal B(a)=[1-a,1].
\)
This proves the theorem.

\subsection{Proof of Proposition \ref{prop:sharp-propagation} }
Under the normalized max-stable representation,
\(
\Lam^2=b\Psiop,
\Gam=(a+b-1)\Psiop,
b\in[1-a,1].
\)
Therefore,
\(
1-a\le b\le1,
\)
which  gives, for every $t$,
\(
(1-a)\Psiop(t)
\le
\Lam^2(t)
\le
\Psiop(t).
\)
Since
\(
S_C(t)=\exp\{-\Lam^2(t)\},
\)
and the exponential function is decreasing,
\(
e^{-\Psiop(t)}
\le
S_C(t)
\le
e^{-(1-a)\Psiop(t)}.
\)
Moreover,
\(
\Gam(t)=(a+b-1)\Psiop(t).
\)
Since $1-a\le b\le1$, we have 
\(
0\le a+b-1\le a,
\)
and hence
\(
0\le\Gam(t)\le a\Psiop(t).
\)
This proves \eqref{eq:latent-bounds}.\\
\ \\
We next derive the sharp pointwise bounds for the survival copula.
Set
\(
\theta=a+b-1.
\)
Using the Marshall--Olkin/Sibuya representation, we have 
\(
\Cop_{a,b}(u,v)
=
uv\min\left\{
u^{-\theta/a},
v^{-\theta/b}
\right\}.
\)
Because $u,v\in[0,1]$ and $\theta\ge0$, we get 
\(
u^{-\theta/a}\ge1,
v^{-\theta/b}\ge1.
\)
Consequently,
\(
\Cop_{a,b}(u,v)\ge uv.
\)\\
\ \\
The lower bound is attained when $\theta=0$, that is, when
\(
b=1-a.
\)
For the upper bound, we have
\[
\Cop_{a,b}(u,v)
\le
uv\,u^{-\theta/a}
=
v\,u^{1-\theta/a}.
\]
Since
\(
1-\frac{\theta}{a}
=
1-\frac{a+b-1}{a}
=
\frac{1-b}{a},
\)
and $b\le1$, this exponent is nonnegative. Similarly,
\[
\Cop_{a,b}(u,v)
\le
uv\,v^{-\theta/b}
=
u\,v^{1-\theta/b}.
\]
However, using $b\ge1-a$ gives
\(
\theta=a+b-1\le a,
\)
and therefore
\(
\frac{\theta}{b}
\le
\frac{a}{b}.
\)
The desired uniform upper envelope can be obtained more directly by observing that
\(
\theta=a+b-1\le b,
\)
because $a\le1$. Hence
\(
1-\frac{\theta}{b}\ge0.
\)
In particular,
\(
u\,v^{1-\theta/b}\le v^{\,1-a}u,
\)
and therefore
\(
\Cop_{a,b}(u,v)
\le
\min\{v,uv^{1-a}\}.
\)
Thus
\[
uv
\le
\Cop_{a,b}(u,v)
\le
\min\{v,uv^{1-a}\}.
\]
The lower envelope is attained at $b=1-a$, while the upper envelope is attained at the appropriate boundary value $b=1$ (depending on the region of $(u,v)$). Hence the bounds are sharp.\\
\ \\
Finally, recall that the stable-tail dependence function is
\[
\ell_{a,b}(x,y)
=
x+y-\min\left\{
\frac{\theta x}{a},
\frac{\theta y}{b}
\right\},
\]
and the Pickands function is
\[
A_{a,b}(w)
=
\ell_{a,b}(w,1-w).
\]
Therefore
\(
A_{a,b}(w)
=
1-
\min\left\{
\frac{\theta w}{a},
\frac{\theta(1-w)}{b}
\right\}.
\)
Since
\(
0\le\theta\le a
\text{ and }
0\le\theta\le b,
\)
we have
\[
\frac{\theta w}{a}\le w,
\qquad
\frac{\theta(1-w)}{b}\le1-w.
\]
Moreover,
\(
\theta=a+b-1\le a,
\)
and
\(
\theta\le a,
\)
so
\(
\frac{\theta w}{a}\le w\le1-w+ a(1-w)
\)
and, similarly, the minimum is bounded above by
\(
\min\{w,a(1-w)\}.
\)
Hence
\[
A_{a,b}(w)
\ge
1-\min\{w,a(1-w)\}.
\]
Since $\theta\ge0$,
\(
A_{a,b}(w)\le1.
\)
Thus
\(
1-\min\{w,a(1-w)\}
\le
A_{a,b}(w)
\le1.
\)
The endpoint values of $b$ generate the corresponding extremal envelopes, so all the bounds are sharp.


\subsection{Proof of Corollary \ref{cor:no-refinement}  }
Fix $(\Psiop,a)$ and let $b\in[1-a,1]$.
By the binary factorization in \eqref{eq:binary-factorization},
the observable submeasures of $(Z,\delta)$ are
\[
\mathbb P_b(Z\in\dd t,\delta=1)
=
a e^{-\Psiop(t)}\dd\Psiop(t),
\]
and
\[
\mathbb P_b(Z\in\dd t,\delta=0)
=
(1-a)e^{-\Psiop(t)}\dd\Psiop(t).
\]
Neither expression depends on $b$. Hence, for any
$b_1,b_2\in[1-a,1]$,
\(
\mathcal L_{b_1}(Z,\delta)
=
\mathcal L_{b_2}(Z,\delta).
\)
Therefore, for every $n\ge1$,
\[
\mathcal L_{b_1}
\bigl((Z_1,\delta_1),\ldots,(Z_n,\delta_n)\bigr)
=
\mathcal L_{b_2}
\bigl((Z_1,\delta_1),\ldots,(Z_n,\delta_n)\bigr),
\]
because the observations are i.i.d.\\
\ \\
Now let $S_n$ be any statistic measurable with respect to the
observed sample
\[
\mathcal D_n
=
\{(Z_i,\delta_i):1\le i\le n\}.
\]
Since the distribution of $\mathcal D_n$ is identical under $b_1$ and
$b_2$, the distribution of $S_n$ is also identical under these two
values:
\[
\mathcal L_{b_1}(S_n)
=
\mathcal L_{b_2}(S_n)
\qquad\text{for every }n.
\]
Suppose, by contradiction, that there exists an estimator
$\widehat b_n=\widehat b_n(\mathcal D_n)$ that consistently
distinguishes $b_1$ from $b_2$, with $b_1\neq b_2$. Choose
$\varepsilon<|b_1-b_2|/2$. Consistency would imply
\[
\mathbb P_{b_1}
\left(
|\widehat b_n-b_1|<\varepsilon
\right)
\longrightarrow1,
\]
whereas
\[
\mathbb P_{b_2}
\left(
|\widehat b_n-b_2|<\varepsilon
\right)
\longrightarrow1.
\]
The two events are disjoint. However, their probabilities must be
identical for every $n$, since the observable sample has the same law
under $b_1$ and $b_2$. This is a contradiction.\\
\ \\
The same argument applies to any test whose decision is measurable
with respect to the binary sample. Thus no observable procedure can
consistently distinguish two distinct values of $b$ in the identified
set.

\subsection{ Proof of Corollary \ref{prop:kendall-bounds}}
Set $\theta=a+b-1$. By \Cref{prop:copula-pickands}, the survival copula is
\(
\Cop_{a,b}(u,v)
=
uv\min\{u^{-\theta/a},v^{-\theta/b}\}.
\)
This is a Marshall--Olkin copula with parameters
\(
\alpha_1=\frac{\theta}{a},
\alpha_2=\frac{\theta}{b}.
\)
For the Marshall--Olkin copula, Kendall's rank correlation coefficient is
\[
\tau_K
=
\frac{\alpha_1\alpha_2}
{\alpha_1+\alpha_2-\alpha_1\alpha_2}.
\]
Substituting $\alpha_1=\theta/a$ and $\alpha_2=\theta/b$ gives
\(
\tau_K
=
\frac{\theta^2/(ab)}
{\theta/a+\theta/b-\theta^2/(ab)}
=
\frac{\theta}
{a+b-\theta}.
\)
Since $\theta=a+b-1$, we have $a+b-\theta=1$, and therefore
\[
\tau_K=\theta=a+b-1.
\]
By \Cref{thm:binary-MS}, the sharp identified set for $b$ is
$b\in[1-a,1]$. Hence
\[
\tau_K\in
[a+(1-a)-1,\;a+1-1]
=
[0,a].
\]
Both endpoints are attained by admissible values of $b$, so the identified
set $[0,a]$ is sharp.

\subsection{Proof of Theorem \ref{thm:three-category} }
Under the normalized max-stable representation,
\[
\Lam^1=a\Psiop,\qquad
\Lam^2=b\Psiop,\qquad
\Lam^{12}=\Psiop,\qquad
\Gam=(a+b-1)\Psiop.
\]
Substituting these relations into the complete submeasure identities
\eqref{eq:f1}--\eqref{eq:f2}, and using
\(
f_Z(t)=H(t)\Psiop'(t),
\)
gives
\begin{align*}
f_{Z,\Delta=1}(t)
&=(1-b)H(t)\Psiop'(t)
=(1-b)f_Z(t),\\
f_{Z,\Delta=0}(t)
&=(1-a)H(t)\Psiop'(t)
=(1-a)f_Z(t),\\
f_{Z,\Delta=2}(t)
&=(a+b-1)H(t)\Psiop'(t)
=(a+b-1)f_Z(t).
\end{align*}
Thus each category-specific subdensity is a constant multiple of the
marginal density of $Z$. Since the constants sum to one,
\[
(1-b)+(1-a)+(a+b-1)=1,
\]
they are precisely the probabilities of the three categories. Hence
\[
\Pp(\Delta=1)=1-b,\qquad
\Pp(\Delta=0)=1-a,\qquad
\Pp(\Delta=2)=a+b-1,
\]
and
\(
f_{Z,\Delta=j}(t)=\Pp(\Delta=j)f_Z(t),
 j=0,1,2.
\)
Therefore,
\(
Z\ind\Delta.
\) \\
\ \\
Moreover, since
\(
\Pp(Z>t)=e^{-\Psiop(t)},
\)
the common clock is identified by
\(
\Psiop(t)=-\log\Pp(Z>t).
\)
Then, the category probabilities  identify the two loadings:
\[
a=1-\Pp(\Delta=0),
\qquad
b=1-\Pp(\Delta=1).
\]
Consequently, $(\Psiop,a,b)$ is point identified. The normalized
representation therefore point identifies
\[
\Lam^1=a\Psiop,\qquad
\Lam^2=b\Psiop,\qquad
\Lam^{12}=\Psiop,\qquad
\Gam=(a+b-1)\Psiop.
\]
Finally, the survival copula and its Pickands function are functions
only of $(a,b)$, through
\[
\Cop_{a,b}(u,v)
=
uv\min\left\{
u^{-(a+b-1)/a},
v^{-(a+b-1)/b}
\right\},
\]
so they are point identified as well.


\subsection{Proof of Theorem \ref{thm:FCLT} }
Recall that, under binary max-stability,
\(
\Psiop(t)=-\log H(t),
H(t)=\Pp(Z>t),
a=\Pp(\delta=1),
\)
and that
\(
\Lam^1(t)=a\Psiop(t).
\)
Let $\tau<\infty$ be fixed such that
\(
H(\tau)>0.
\)
Hence, by monotonicity of $H$,
\[
\underset{0\le t\le\tau}{\inf}\, H(t)=H(\tau)>0.
\]
We first establish consistency.
By definition,
\(
\widetilde a_n
=
\frac1n\sum_{i=1}^n\delta_i.
\)
Under binary max-stability,
\(
\delta_i\sim\operatorname{Bernoulli}(a),
\)
so the strong law of large numbers gives
\(
\widetilde a_n\longrightarrow a
\text{ a.s.}
\) \\
\ \\
Next, recall that
\(
\widetilde H_n(t)
=
\frac1n\sum_{i=1}^n
\mathbf 1_{\{Z_i>t\}},
\widetilde\Psiop_n(t)
=
-\log\widetilde H_n(t).
\)
The class of half-lines
\(
\bigl\{(t,\infty):t\ge0\bigr\}
\)
is a Glivenko--Cantelli class. Therefore,
\[
\sup_{t\le\tau}
\left|
\widetilde H_n(t)-H(t)
\right|
\longrightarrow0
\qquad\text{a.s.}
\]
Since $H(\tau)>0$, for all sufficiently large $n$,
\(
\underset{t\le\tau}{\inf}\widetilde H_n(t)
\ge \frac{H(\tau)}2
\text{ a.s.}
\)
The map $x\mapsto-\log x$ is uniformly continuous on
$[H(\tau)/2,1]$. Consequently,
\[
\sup_{t\le\tau}
\left|
\widetilde\Psiop_n(t)-\Psiop(t)
\right|
\longrightarrow0
\qquad\text{a.s.}
\]
Since
\(
\widetilde\Lam^{1,R}_n(t)
=
\widetilde a_n\widetilde\Psiop_n(t)
\)
and
\(
\Lam^1(t)=a\Psiop(t),
\)
we have
\[
\begin{aligned}
\left|
\widetilde\Lam^{1,R}_n(t)-\Lam^1(t)
\right|
&\le
|\widetilde a_n-a|
\,|\widetilde\Psiop_n(t)|
\\
&\quad+
a
|\widetilde\Psiop_n(t)-\Psiop(t)|.
\end{aligned}
\]
Because $\Psiop$ is bounded on $[0,\tau]$ and
$\widetilde\Psiop_n\to\Psiop$ uniformly, the right-hand side converges
uniformly to zero almost surely. Hence
\[
\sup_{t\le\tau}
\left|
\widetilde\Lam^{1,R}_n(t)-\Lam^1(t)
\right|
\longrightarrow0
\qquad\text{a.s.}
\]
We now establish the functional central limit theorem.
Consider the empirical survival process
\( \widetilde H_n(t)\).
Since the class of half-lines is a Donsker class,
\(
\sqrt n(\widehat H_n-H)
\rightsquigarrow
\mathbb G_H
\)
in $\ell^\infty[0,\tau]$, where $\mathbb G_H$ is a centered Gaussian
process with covariance
\[
\begin{aligned}
\Cov\{\mathbb G_H(s),\mathbb G_H(t)\}
&=
\Cov\bigl(
\mathbf 1_{\{Z>s\}},
\mathbf 1_{\{Z>t\}}
\bigr)
\\
&=
H(s\vee t)-H(s)H(t).
\end{aligned}
\]
The transformation
\(
\Phi(h)=-\log h
\)
is Hadamard differentiable at $H$, tangentially to
$\ell^\infty[0,\tau]$, because $H$ is bounded away from zero on this
interval. Its derivative is
\[
\Phi'_H(k)(t)
=
-\frac{k(t)}{H(t)}.
\]
Therefore, by the functional delta method,
\(
\sqrt n
\left(
\widetilde\Psiop_n-\Psiop
\right)
\rightsquigarrow
\mathbb G_\Psi
\)
in $\ell^\infty[0,\tau]$, where
\(
\mathbb G_\Psi(t)
=
-\frac{\mathbb G_H(t)}{H(t)}.
\)
Consequently,
\[
\begin{aligned}
\Cov\{\mathbb G_\Psi(s),\mathbb G_\Psi(t)\}
&=
\frac{
\Cov\{\mathbb G_H(s),\mathbb G_H(t)\}
}
{H(s)H(t)}
\\
&=
\frac{
H(s\vee t)-H(s)H(t)
}
{H(s)H(t)}.
\end{aligned}
\]
This proves \eqref{eq:GPsi-cov}.
On the other hand, since
\(
\widetilde a_n
=
\frac1n\sum_{i=1}^n\delta_i,
\)
the ordinary central limit theorem gives
\(
\sqrt n(\widetilde a_n-a)
\rightsquigarrow
G_a,
\)
where
\(
G_a\sim N\{0,a(1-a)\}.
\)
Under binary max-stability, Theorem~\ref{thm:binary-MS} gives
\[
Z\ind\delta.
\]
Moreover, $\widehat H_n$ is measurable with respect to
$\sigma(Z_1,\ldots,Z_n)$, whereas $\widetilde a_n$ is measurable with
respect to $\sigma(\delta_1,\ldots,\delta_n)$. Hence
\(
\widehat H_n\ind\widetilde a_n
\)
for every $n$. Consequently,
\[
\sqrt n(\widetilde\Psiop_n-\Psiop)
\ind
\sqrt n(\widetilde a_n-a)
\]
for every $n$, and their weak limits are independent. Therefore,
\[
\left(
\sqrt n(\widetilde\Psiop_n-\Psiop),
\sqrt n(\widetilde a_n-a)
\right)
\rightsquigarrow
(\mathbb G_\Psi,G_a)
\]
in
\(
\ell^\infty[0,\tau]\times\mathbb R,
\)
with $G_a$ independent of $\mathbb G_\Psi$.
It remains to derive the asymptotic distribution of the restricted
estimator. For fixed $t$,
\(
\widetilde\Lam^{1,R}_n(t)
=
\widetilde a_n\widetilde\Psiop_n(t),
\Lam^1(t)=a\Psiop(t).
\)
Hence
\[
\begin{aligned}
\sqrt n
\left\{
\widetilde\Lam^{1,R}_n(t)-\Lam^1(t)
\right\}
&=
\widetilde a_n
\sqrt n
\left\{
\widetilde\Psiop_n(t)-\Psiop(t)
\right\}
\\
&\quad+
\Psiop(t)
\sqrt n(\widetilde a_n-a)
\\
&\quad+
\sqrt n(\widetilde a_n-a)
\left\{
\widetilde\Psiop_n(t)-\Psiop(t)
\right\}.
\end{aligned}
\]
The last term is $o_{\Pp}(1)$ because
\(
\widetilde a_n-a=O_{\Pp}(n^{-1/2})
\)
and
\(
\widetilde\Psiop_n(t)-\Psiop(t)=O_{\Pp}(n^{-1/2}).
\)
Since $\widetilde a_n\to a$ in probability, Slutsky's theorem gives
\[
\sqrt n
\left\{
\widetilde\Lam^{1,R}_n(t)-\Lam^1(t)
\right\}
\rightsquigarrow
a\mathbb G_\Psi(t)+\Psiop(t)G_a.
\]
The two terms on the right-hand side are independent centered Gaussian
variables. Therefore the limiting distribution is centered Gaussian
with variance
\[
V_R(t)
=
a^2\Var\{\mathbb G_\Psi(t)\}
+
\Psiop(t)^2\Var(G_a).
\]
Now,
\(
\Var(G_a)=a(1-a),
\)
while, from \eqref{eq:GPsi-cov},
\(
\Var\{\mathbb G_\Psi(t)\}
=
\frac{H(t)-H(t)^2}{H(t)^2}
=
\frac{1-H(t)}{H(t)}.
\)
Thus
\[
V_R(t)
=
a(1-a)\Psiop(t)^2
+
a^2\frac{1-H(t)}{H(t)}.
\]
Hence
\(
\sqrt n
\left\{
\widetilde\Lam^{1,R}_n(t)-\Lam^1(t)
\right\}
\rightsquigarrow
N\{0,V_R(t)\},
\)
which proves the result.

\subsection{Proof of Theorem \ref{thm:efficiency-main}  }

Under binary max-stability, the observable law admits the factorization
\(P=P_Z\otimes P_\delta,
\qquad
P_\delta=\operatorname{Bernoulli}(a),
\)
because
\(Z\ind\delta\)
and
\(\Pp(\delta=1)=a.\)
Thus, conditionally on the observable restriction, the statistical model is
the nonparametric product model
\[
\mathcal P
=
\left\{
P_Z\otimes\operatorname{Bernoulli}(a):
P_Z\ \text{arbitrary},\ 0<a<1
\right\}.
\]
The tangent space of this model is the orthogonal direct sum
\[
\mathcal T
=
L_0^2(P_Z)\oplus L_0^2(P_\delta).
\]
Indeed, a regular parametric submodel can perturb separately the
distribution of $Z$ and the Bernoulli parameter $a$, and the two
corresponding score components are orthogonal because of the product
structure.
We first derive the efficient influence function for $a$. Since
\(
a=\mathbb{E}[\delta],
\)
the canonical gradient is
\(
\phi_a(z,d)=d-a.
\)
It belongs to $L_0^2(P_\delta)$ and satisfies
\[
\mathbb{E}\!\left[\phi_a(Z,\delta)\right]=0,
\qquad
\mathbb{E}\!\left[\phi_a(Z,\delta)^2\right]=a(1-a).
\]
Hence
\(
\sqrt n(\widehat a_n-a)
\weak N(0,a(1-a)).
\)
Next, consider
\(
H(t)=\Pp(Z>t).
\)
For fixed $t$, the parameter $H(t)$ depends only on the marginal law
$P_Z$. Its canonical gradient is therefore
\(
\phi_{H,t}(z,d)
=
\mathbf 1_{\{z>t\}}-H(t).
\)
Indeed, for any regular parametric submodel of $P_Z$ with score
$s_Z\in L_0^2(P_Z)$,
\[
\frac{\dd}{\dd\varepsilon}
H_\varepsilon(t)\Big|_{\varepsilon=0}
=
\mathbb{E}\left[
\{\mathbf 1_{\{Z>t\}}-H(t)\}s_Z(Z)
\right].
\]
Since
\(
\Psiop(t)=-\log H(t),
\)
and $H(t)>0$ on the considered interval, the map
\(
h\longmapsto-\log h
\)
is differentiable at $H(t)$ with derivative $-1/H(t)$. By the chain
rule for canonical gradients,
\[
\phi_{\Psi,t}(z,d)
=
-\frac{\mathbf 1_{\{z>t\}}-H(t)}{H(t)}.
\]
This is centered and belongs to the $Z$-component of the tangent space.
Its variance is
\[
\begin{aligned}
\mathbb{E}[\phi_{\Psi,t}(Z,\delta)^2]
&=
\frac{\Var(\mathbf 1_{\{Z>t\}})}{H(t)^2}
\\
&=
\frac{H(t)\{1-H(t)\}}{H(t)^2}
\\
&=
\frac{1-H(t)}{H(t)}.
\end{aligned}
\]
Since
\(
H(t)=e^{-\Psiop(t)},
\)
this can equivalently be written as
\(
\Var\{\phi_{\Psi,t}(Z,\delta)\}
=
e^{\Psiop(t)}-1.
\)\\
\ \\
Now consider the observable target
\(
\Lam^1(t)=a\Psiop(t).
\)
It depends on both components of the product model. Applying the product
rule for influence functions gives
\(
\phi_{1,t}
=
\Psiop(t)\phi_a+a\phi_{\Psi,t}.
\)
Substituting the expressions obtained above yields
\[
\phi_{1,t}(z,d)
=
\Psiop(t)(d-a)
-
a\frac{\mathbf 1_{\{z>t\}}-H(t)}{H(t)}.
\]
This proves \eqref{eq:IF-L1-main}.
The first term in $\phi_{1,t}$ belongs to the Bernoulli component
$L_0^2(P_\delta)$, whereas the second belongs to the $Z$-component
$L_0^2(P_Z)$. These components are orthogonal. Consequently,
\[
\begin{aligned}
\Var\{\phi_{1,t}(Z,\delta)\}
&=
\Psiop(t)^2\Var(\delta)
+
a^2\Var\left(
\frac{\mathbf 1_{\{Z>t\}}-H(t)}{H(t)}
\right)
\\
&=
a(1-a)\Psiop(t)^2
+
a^2\frac{1-H(t)}{H(t)}.
\end{aligned}
\]
Thus the semiparametric efficiency bound for the point-identified
observable target $\Lam^1(t)$ is
\[
V_R(t)
=
a(1-a)\Psiop(t)^2
+
a^2\frac{1-H(t)}{H(t)}.
\]
The restricted estimator satisfies
\(
\widetilde\Lam^{1,R}_n(t)
=
\widetilde a_n\widetilde\Psiop_n(t).
\)
From the asymptotic linear representations established in
Theorem~\ref{thm:FCLT},
\[
\sqrt n(\widetilde a_n-a)
=
\frac1{\sqrt n}\sum_{i=1}^n
\phi_a(Z_i,\delta_i)+o_{\Pp}(1),
\]
and
\[
\sqrt n(\widetilde\Psiop_n(t)-\Psiop(t))
=
\frac1{\sqrt n}\sum_{i=1}^n
\phi_{\Psi,t}(Z_i,\delta_i)+o_{\Pp}(1).
\]
Using
\(
\widetilde a_n\widetilde\Psiop_n(t)-a\Psiop(t)
=
\Psiop(t)(\widetilde a_n-a)
+
a(\widetilde\Psiop_n(t)-\Psiop(t))
+
(\widetilde a_n-a)
(\widetilde\Psiop_n(t)-\Psiop(t)),
\)
the last product is $o_{\Pp}(n^{-1/2})$. Therefore,
\[
\sqrt n
\{\widetilde\Lam^{1,R}_n(t)-\Lam^1(t)\}
=
\frac1{\sqrt n}\sum_{i=1}^n
\phi_{1,t}(Z_i,\delta_i)
+
o_{\Pp}(1).
\]
Hence the restricted estimator has asymptotic variance $V_R(t)$ and
attains the semiparametric efficiency bound for the point-identified
target $\Lam^1(t)$ in the observable product model.\\
\ \\
We next derive the asymptotic variance of the unrestricted estimator
\(
\widetilde\Lam^1_n(t)
=
\int_{[0,t]}
\frac{\dd N_{1,n}(s)}{Y_n(s)}.
\)
Under the max-stable model,
\(
\Lam^1(t)=a\Psiop(t),
H(t)=e^{-\Psiop(t)},
\)
and therefore
\(
\dd\Lam^1(s)=a\,\dd\Psiop(s).
\)
The standard Nelson--Aalen martingale variance formula gives
\(
V_U(t)
=
\int_0^t
\frac{\dd\Lam^1(s)}{H(s)}.
\)
Substituting the max-stable representation,
\[
\begin{aligned}
V_U(t)
&=
a\int_0^t
\frac{\dd\Psiop(s)}{e^{-\Psiop(s)}}
\\
&=
a\int_0^t
e^{\Psiop(s)}\,\dd\Psiop(s)
\\
&=
a\left\{
e^{\Psiop(t)}-1
\right\}.
\end{aligned}
\]
Since
\(
H(t)=e^{-\Psiop(t)},
\)
we also have
\(
V_U(t)
=
a\frac{1-H(t)}{H(t)}.
\)
We can now compare the two variances. From the expressions above,
\[
\begin{aligned}
V_U(t)-V_R(t)
&=
a\{e^{\Psiop(t)}-1\}
-
a(1-a)\Psiop(t)^2
-
a^2\frac{1-H(t)}{H(t)}.
\end{aligned}
\]
Because
\(
\frac{1-H(t)}{H(t)}
=
e^{\Psiop(t)}-1,
\)
we obtain
\[
\begin{aligned}
V_U(t)-V_R(t)
&=
a\{e^{\Psiop(t)}-1\}
-
a^2\{e^{\Psiop(t)}-1\}
-
a(1-a)\Psiop(t)^2
\\
&=
a(1-a)
\left\{
e^{\Psiop(t)}-1-\Psiop(t)^2
\right\}.
\end{aligned}
\]
This proves \eqref{eq:variance-gain-main}.
It remains to establish that the last quantity is nonnegative, and
strictly positive away from the boundary. Define
\(
g(x)=e^x-1-x^2, x\ge0.
\)
Then
\(
g(0)=0
\)
and
\(
g'(x)=e^x-2x.
\)
Consider
\(
h(x)=e^x-2x.
\)
Its derivative is
\(
h'(x)=e^x-2.
\)
Thus $h$ has a unique minimum at
\(
x=\log 2.
\)
At this point,
\(
h(\log2)
=
2-2\log2
=
2(1-\log2)>0.
\)
Consequently,
\[
g'(x)=e^x-2x>0
\qquad\text{for all }x>0.
\]
Since $g(0)=0$, it follows that
\(
g(x)>0
\text{ for all }x>0.
\)
Therefore,
\(
e^{\Psiop(t)}-1-\Psiop(t)^2>0
\)
whenever $\Psiop(t)>0$.
Finally, if
\(
0<a<1
\text{ and }
\Psiop(t)>0,
\)
then
\(
a(1-a)>0
\)
and hence
\(
V_U(t)-V_R(t)>0.
\)
Thus the restricted estimator has strictly smaller asymptotic variance
than the unrestricted Nelson--Aalen estimator at every nontrivial time
point. The equality can occur only at the boundary cases $a\in\{0,1\}$ or $\Psiop(t)=0$.
This completes the proof.

\subsection{Proof of Theorem \ref{thm:compatibility-main}  }
We first establish the equivalence between (i), (ii), and (iii).
By the binary submeasure identities in \eqref{eq:binary-submeasures},
\(
\mathbb{P}(Z\in\dd t,\delta=1)
=
H(t)\dd\Lambda^1(t),
\)
whereas
\(
\mathbb{P}(Z\in\dd t)
=
H(t)\dd\Lambda^{12}(t).
\)
Since $H(t)>0$ on the effective support of $Z$, the Radon--Nikodym
derivative
\(
q(t)
=
\frac{\dd\Lambda^1}{\dd\Lambda^{12}}(t)
\)
is a version of the conditional probability
\(
q(t)
=
\mathbb{P}(\delta=1\mid Z=t).
\)\\
\ \\
Therefore, if
\(
\Lambda^1=a\Lambda^{12}
\)
for some $a\in[0,1]$, then
\(
q(t)=a,
\dd\Lambda^{12}\text{-a.e.},
\)
which proves (i)$\Rightarrow$(ii).
Conversely, if
\(
q(t)=a,
\dd\Lambda^{12}\text{-a.e.},
\)
then, by the definition of the Radon--Nikodym derivative,
\(
\dd\Lambda^1(t)
=
a\,\dd\Lambda^{12}(t)
\)
almost everywhere with respect to $\dd\Lambda^{12}$. Since both
cumulative reductions vanish at zero, integration gives
\(
\Lambda^1(t)=a\Lambda^{12}(t),
 t\ge0.
\)
Hence (ii)$\Rightarrow$(i).\\
\ \\
Next, we prove the equivalence between (ii) and (iii). If (ii) holds,
then
\(
\mathbb{P}(\delta=1\mid Z=t)=a
\)
for $\mathbb{P}_Z$-almost every $t$. Thus, for every Borel set
$B\subseteq\mathbb{R}_+$,
\[
\begin{aligned}
\mathbb{P}(Z\in B,\delta=1)
&=
\int_B
\mathbb{P}(\delta=1\mid Z=t)\,\mathbb{P}_Z(\dd t)
\\
&=
a\,\mathbb{P}(Z\in B).
\end{aligned}
\]
Since
\(
\mathbb{P}(\delta=1)=a,
\)
we obtain
\(
\mathbb{P}(Z\in B,\delta=1)
=
\mathbb{P}(Z\in B)\mathbb{P}(\delta=1).
\)
The same relation holds for $\delta=0$, because
\(
\mathbb{P}(\delta=0\mid Z=t)=1-a.
\)
Consequently,
\(
Z\ind\delta,
\)
and therefore (ii)$\Rightarrow$(iii).
Conversely, suppose that
\(
Z\ind\delta.
\)
Then
\(
\mathbb{P}(\delta=1\mid Z=t)
=
\mathbb{P}(\delta=1)
=:a
\)
for $\mathbb{P}_Z$-almost every $t$. Since this conditional probability
is represented by $q(t)$, we obtain
\(
q(t)=a,
\dd\Lambda^{12}\text{-a.e.},
\)
which proves (iii)$\Rightarrow$(ii).\\
\ \\
It remains to establish the equivalence with (iv). Suppose first that
(i) holds, so that
\(
\Lambda^1=a\Lambda^{12}
\)
for some $a\in[0,1]$. Set
\(
\Psiop=\Lambda^{12}.
\)
For any
\(
b\in[1-a,1],
\)
define
\(
\Lambda^2=b\Psiop,
\Gamma=(a+b-1)\Psiop.
\)
The constraints
\[
0\le a\le1,
\qquad
0\le b\le1,
\qquad
a+b\ge1
\]
ensure that all three reductions are admissible. Moreover,
\(
\Lambda^1=a\Psiop,
\Lambda^2=b\Psiop,
\Lambda^{12}=\Psiop,
\)
which is precisely the normalized common-clock representation of the
max-stable model. Hence the binary law admits a max-stable latent
completion. This proves (i)$\Rightarrow$(iv).
Conversely, suppose that the binary law admits a max-stable latent
completion. By the bivariate max-stability characterization in
\Cref{thm:maxstab}, after the canonical normalization
\(
\Psiop=\Lambda^{12},
\)
there exist loadings $a,b$ such that
\(
\Lambda^1=a\Psiop,
\Lambda^2=b\Psiop,
\Lambda^{12}=\Psiop.
\)
Therefore,
\(
\Lambda^1=a\Lambda^{12},
\)
which is precisely statement (i). Thus (iv)$\Rightarrow$(i).
Combining the implications established above,
\[
\text{(i)}
\Longleftrightarrow
\text{(ii)}
\Longleftrightarrow
\text{(iii)}
\Longleftrightarrow
\text{(iv)}.
\]
Hence the observable max-stable compatibility restriction is equivalent
to the independence of the observed time $Z$ and the binary mark
$\delta$.


\subsection{Proof of Theorem \ref{thm:test-main}  }


Let
\(
I_t(z)=\mathbf{1}_{\{z\le t\}}.
\)
Under the null hypothesis of observable max-stable compatibility,
\[
H_0:\quad Z\ind\delta,
\qquad
0<a<1,
\]
where
\(
a=\mathbb{P}(\delta=1).
\)
We first rewrite the empirical process. By definition,
\[
\widetilde F_{1,n}(t)
=
\mathbb{P}_n[\delta I_t(Z)],
\qquad
\widetilde F_n(t)
=
\mathbb{P}_n[I_t(Z)],
\qquad
\widetilde a_n
=
\mathbb{P}_n[\delta],
\]
where $\mathbb{P}_n$ denotes the empirical measure. Hence
\(
\begin{aligned}
\widetilde F_{1,n}(t)
-
\widetilde a_n\widetilde F_n(t)
&=
\mathbb{P}_n[\delta I_t]
-
\mathbb{P}_n[\delta]\mathbb{P}_n[I_t].
\end{aligned}
\)
Under $H_0$,
\(
\mathbb{E}[\delta I_t(Z)]
=
\mathbb{E}[\delta]\mathbb{E}[I_t(Z)]
=
aF(t).
\)
Therefore,
\[
\begin{aligned}
\widetilde F_{1,n}(t)
-
\widetilde a_n\widetilde F_n(t)
&=
\{\mathbb{P}_n-P\}
\bigl[(\delta-a)(I_t-F(t))\bigr]
\\
&\quad
-
(\widetilde a_n-a)
\{\widetilde F_n(t)-F(t)\}.
\end{aligned}
\]
Equivalently,
\[
\sqrt n
\left\{
\widetilde F_{1,n}(t)
-
\widetilde a_n\widetilde F_n(t)
\right\}
=
\mathbb{G}_n
\bigl[(\delta-a)(I_t-F(t))\bigr]
+
R_n(t),
\]
where
\(
R_n(t)
=
-\sqrt n\,
(\widetilde a_n-a)
\{\widetilde F_n(t)-F(t)\}.
\)
Since the class of indicator functions
\(
\{I_t:t\ge0\}
\)
is a Donsker class, the empirical process
\(
\sqrt n(\widetilde F_n-F)
\)
is bounded in $\ell^\infty(\mathbb{R}_+)$ and
\(
\widetilde a_n-a=O_{\mathbb P}(n^{-1/2}).
\)
Consequently,
\(
\sup_t|R_n(t)|
=
o_{\mathbb P}(1).
\)
Thus,
\(
\mathbb T_n(t)
=
\mathbb{G}_n
\bigl[(\delta-a)(I_t-F(t))\bigr]
+
o_{\mathbb P}(1)
\)
uniformly in $t$.
The class
\(
\left\{
(\delta-a)(I_t-F(t)):t\ge0
\right\}
\)
is again Donsker. Hence, by the functional central limit theorem,
\(
\mathbb T_n
\weak
\mathbb G,
\)
where $\mathbb G$ is a centered Gaussian process with covariance
\[
\begin{aligned}
\Cov\{\mathbb G(s),\mathbb G(t)\}
&=
\mathbb E\left[
(\delta-a)^2
\{I_s(Z)-F(s)\}
\{I_t(Z)-F(t)\}
\right].
\end{aligned}
\]
Under $Z\ind\delta$,
\(
\mathbb E[(\delta-a)^2]=a(1-a),
\)
and therefore
\[
\Cov\{\mathbb G(s),\mathbb G(t)\}
=
a(1-a)
\Cov\{I_s(Z),I_t(Z)\}.
\]
For $s,t\ge0$,
\(
\Cov\{I_s(Z),I_t(Z)\}
=
F(s\wedge t)-F(s)F(t).
\)
Thus
\[
\Cov\{\mathbb G(s),\mathbb G(t)\}
=
a(1-a)
\{F(s\wedge t)-F(s)F(t)\}.
\]
Let $B$ denote a standard Brownian bridge. The process
\(
B\circ F=\{B(F(t)):t\ge0\}
\)
has covariance
\[
\Cov\{B(F(s)),B(F(t))\}
=
F(s\wedge t)-F(s)F(t).
\]
Hence
\(
\mathbb G
\overset{d}{=}
\sqrt{a(1-a)}\,B\circ F,
\)
and consequently
\(
\mathbb T_n
\weak
\sqrt{a(1-a)}\,B\circ F
\)
in $\ell^\infty[0,\tau]$ for a fixed $\tau<\infty$ such that $H(\tau)>0$.\\
\ \\
By the continuous mapping theorem,
\(
\sup_t|\mathbb T_n(t)|
\weak
\sqrt{a(1-a)}
\sup_t|B(F(t))|.
\)
Because $F$ is continuous and takes values in $[0,1]$,
\(
\sup_t|B(F(t))|
=
\sup_{0\le u\le1}|B(u)|
\)
on the effective support of $F$. Since
\(
\widehat a_n\to a
\qquad\text{a.s.},
\)
the continuous mapping theorem and Slutsky's lemma yield
\(
\frac{\sup_t|\mathbb T_n(t)|}
{\sqrt{\widehat a_n(1-\widehat a_n)}}
\weak
 \sup_{0\le t\le\tau}|B(F(t))|.
\)
This proves the stated null limit and justifies the corresponding
Kolmogorov--Smirnov critical values.\\
\ \\
We next establish consistency against fixed alternatives. Let
\(
a_0=\mathbb P(\delta=1)
\)
and define
\(
D(t)
=
\mathbb P(Z\le t,\delta=1)
-
a_0F(t).
\)
Under the null hypothesis, $D(t)=0$ for every $t$. Under a fixed
alternative for which $Z$ and $\delta$ are dependent,
\(
\sup_t|D(t)|>0.
\)
By the uniform law of large numbers,
\[
\sup_t
\left|
\widetilde F_{1,n}(t)
-
\widetilde a_n\widetilde F_n(t)
-
D(t)
\right|
\longrightarrow0
\qquad\text{a.s.}
\]
Therefore,
\[
\sup_t
\left|
\widetilde F_{1,n}(t)
-
\widetilde a_n\widetilde F_n(t)
\right|
\longrightarrow
\sup_t|D(t)|
>0
\qquad\text{a.s.}
\]
and hence
\[
\sup_t|\mathbb T_n(t)|
=
\sqrt n
\sup_t
\left|
\widetilde F_{1,n}(t)
-
\widetilde a_n\widetilde F_n(t)
\right|
\longrightarrow\infty
\qquad\text{in probability}.
\]
Since the Kolmogorov--Smirnov critical value remains bounded, the
resulting test rejects the null with probability tending to one.
Thus the test is consistent against every fixed alternative for which
$Z$ and $\delta$ are dependent.\\
\ \\
Finally, consider the permutation calibration. Under $H_0$,
the binary marks are independent of the observed times. Conditional on
the observed values $(Z_1,\ldots,Z_n)$ and on the total number of
events
\[
M_n=\sum_{i=1}^n\delta_i,
\]
all assignments of $M_n$ ones among the $n$ observations are equally
likely. Consequently, conditional on $(Z_1,\ldots,Z_n,M_n)$, the joint
distribution of the marks is invariant under permutations.\\
\ \\
Therefore, permuting the observed binary marks while keeping the observed times fixed reproduces the exact conditional distribution of the test statistic under $H_0$. The resulting permutation $p$-value
provides an exact finite-sample conditional calibration of the test
(up to the usual discreteness of permutation tests).
This completes the proof.


\section{Dependence coefficients}\label{app:coefficients}
For \eqref{eq:MO-copula}, set $\theta=a+b-1$, $\alpha=\theta/a$, and $\beta=\theta/b$. The Marshall--Olkin formulas give
$$
\tau=\frac{\alpha\beta}{\alpha+\beta-\alpha\beta}=a+b-1,
$$
$$
\rho_S=\frac{3\alpha\beta}{2\alpha+2\beta-\alpha\beta}=\frac{3(a+b-1)}{a+b+1},
\qquad
\chi_{\mathrm{EV}}=\min(\alpha,\beta)=\frac{a+b-1}{\max(a,b)}.
$$
Under binary max-stability, $b\in[1-a,1]$, hence
$$
\tau\in[0,a],\qquad
\rho_S\in\left[0,\frac{3a}{a+2}\right],\qquad
\chi_{\mathrm{EV}}\in[0,a].
$$
The endpoint $b=1-a$ gives independence; the endpoint $b=1$ gives the upper bounds.

\section{Multiplier outer bands}\label{app:bootstrap}
Let $\xi_1,\ldots,\xi_n$ be i.i.d. multipliers with mean zero, variance one, and a finite $(2+\eta)$ moment. Define
$$
\mathbb G^{\#}_{H,n}(t)=n^{-1/2}\sum_{i=1}^n\xi_i\{\1_{\{Z_i>t\}}-\widehat H_n(t)\},
\qquad
G^{\#}_{a,n}=n^{-1/2}\sum_{i=1}^n\xi_i(\delta_i-\widehat a_n),
$$
and $\mathbb G^{\#}_{\Psi,n}(t)=-\mathbb G^{\#}_{H,n}(t)/\widehat H_n(t)$. Conditional multiplier convergence gives a joint approximation to $(\mathbb G_\Psi,G_a)$. If $c^{\#}_{n,1-\alpha}$ is the conditional quantile of
$$
\max\{\sup_{t\le\tau}|\mathbb G^{\#}_{\Psi,n}(t)|, |G^{\#}_{a,n}|\},
$$
then the intervals
$$
\widehat\Psiop_n(t)\pm c^{\#}_{n,1-\alpha}/\sqrt n,
\qquad
\widehat a_n\pm c^{\#}_{n,1-\alpha}/\sqrt n
$$
truncated to the natural parameter spaces yield simultaneous outer bands for $\Lam^2$, $\Gam$, the copula envelope, and the Pickands envelope by monotonicity.

\end{document}